\documentclass[12pt]{article}
\usepackage{amsfonts}
\usepackage{bbm}
\usepackage{amssymb}
\newdimen\DiagramCellHeight\DiagramCellHeight3em 
\newdimen\DiagramCellWidth\DiagramCellWidth3em 
\newdimen\MapBreadth\MapBreadth.04em 
\newdimen\MapShortFall\MapShortFall.4em 
\newdimen\PileSpacing\PileSpacing1em 
\def\labelstyle{\ifincommdiag\textstyle\else\scriptstyle\fi}
\let\objectstyle\displaystyle
\def\NE{\NorthEast\DiagonalMap{\lah22}{\laf0}{}{\laf0}{}(2,2)}

\def\rTo{\HorizontalMap\empty-\empty-\rhvee}
\def\lTo{\HorizontalMap\lhvee-\empty-\empty}
\def\dTo{\VerticalMap\empty|\empty|\dhvee}

\font\tenln=line10

\mathchardef\lt="313C \mathchardef\gt="313E

\def\rhvee{\mkern-10mu\gt}
\def\lhvee{\lt\mkern-10mu}
\def\dhvee{\vbox\tozpt{\vss\hbox{$\vee$}\kern0pt}}
\def\dhcvee{\vbox\tozpt{\vss\hbox{$\curlyvee$}\kern0pt}}
\def\dhvvee{\vbox\tozpt{\vss\hbox{$\vee$}\kern-.6ex\hbox{$\vee$}\kern0pt}}
\def\uhvvee{\vbox\tozpt{\hbox{$\wedge$}\kern-.6ex\hbox{$\wedge$}\vss}}
\def\twoheaddownarrow{\rlap{$\downarrow$}\raise-.5ex\hbox{$\downarrow$}}
\def\twoheaduparrow{\rlap{$\uparrow$}\raise.5ex\hbox{$\uparrow$}}
\def\rhla{\vbox\tozpt{\vss\hbox\tozpt{\hss\tenln\char'55}\kern\axisheight}}
\def\lhla{\vbox\tozpt{\vss\hbox\tozpt{\tenln\char'33\hss}\kern\axisheight}}
\def\rthooka{\raise.603ex\hbox{$\scriptscriptstyle\subset$}}
\def\lthooka{\raise.603ex\hbox{$\scriptscriptstyle\supset$}}
\def\rthookb{\raise-.022ex\hbox{$\scriptscriptstyle\subset$}}
\def\lthookb{\raise-.022ex\hbox{$\scriptscriptstyle\supset$}}

\def\SEpbk{\rlap{\smash{\kern0.1em \vrule depth 2.67ex height -2.55ex width 0%
.9em \vrule height -0.46ex depth 2.67ex width .05em }}}
\def\SWpbk{\llap{\smash{\vrule height -0.46ex depth 2.67ex width .05em \vrule
depth 2.67ex height -2.55ex width .9em \kern0.1em }}}
\def\NEpbk{\rlap{\smash{\kern0.1em \vrule depth -3.48ex height 3.67ex width 0%
.95em \vrule height 3.67ex depth -1.39ex width .05em }}}
\def\NWpbk{\llap{\smash{\vrule height 3.6ex depth -1.39ex width .05em \vrule
depth -3.48ex height 3.67ex width .95em \kern0.1em }}}

\newcount\cdna\newcount\cdnb\newcount\cdnc\newcount\cdnd\cdna=\catcode`\@%
\catcode`\@=11 \let\then\relax\def\loopa#1\repeat{\def\bodya{#1}\iteratea}%
\def\iteratea{\bodya\let\next\iteratea\else\let\next\relax\fi\next}\def\loopb
#1\repeat{\def\bodyb{#1}\iterateb}\def\iterateb{\bodyb\let\next\iterateb\else
\let\next\relax\fi\next} \def\mapctxterr{\message{commutative diagram: map
context error}}\def\mapclasherr{\message{commutative diagram: clashing maps}}%
\def\ObsDim#1{\expandafter\message{! diagrams Warning: Dimension \string#1 is
obsolete (ignored)}\global\let#1\ObsDimq\ObsDimq}\def\ObsDimq{\dimen@=}\def
\HorizontalMapLength{\ObsDim\HorizontalMapLength}\def\VerticalMapHeight{%
\ObsDim\VerticalMapHeight}\def\VerticalMapDepth{\ObsDim\VerticalMapDepth}\def
\VerticalMapExtraHeight{\ObsDim\VerticalMapExtraHeight}\def
\VerticalMapExtraDepth{\ObsDim\VerticalMapExtraDepth}\def\ObsCount#1{%
\expandafter\message{! diagrams Warning: Count \string#1 is obsolete (ignored%
)}\global\let#1\ObsCountq\ObsCountq}\def\ObsCountq{\count@=}\def
\DiagonalLineSegments{\ObsCount\DiagonalLineSegments}\def\tozpt{to\z@}\def
\sethorizhtdp{\dimen8=\axisheight\dimen9=\MapBreadth\advance\dimen8.5\dimen9%
\advance\dimen9-\dimen8}\def\horizhtdp{height\dimen8 depth\dimen9 }\def
\axisheight{\fontdimen22\the\textfont2 }\countdef\boxc@unt=14

\def\bombparameters{\hsize\z@\rightskip\z@ plus1fil minus\maxdimen
\parfillskip\z@\linepenalty9000 \looseness0 \hfuzz\maxdimen\hbadness10000
\clubpenalty0 \widowpenalty0 \displaywidowpenalty0 \interlinepenalty0
\predisplaypenalty0 \postdisplaypenalty0 \interdisplaylinepenalty0
\interfootnotelinepenalty0 \floatingpenalty0 \brokenpenalty0 \everypar{}%
\leftskip\z@\parskip\z@\parindent\z@\pretolerance10000 \tolerance10000
\hyphenpenalty10000 \exhyphenpenalty10000 \binoppenalty10000 \relpenalty10000
\adjdemerits0 \doublehyphendemerits0 \finalhyphendemerits0 \prevdepth\z@}\def
\startbombverticallist{\hbox{}\penalty1\nointerlineskip}

\def\pushh#1\to#2{\setbox#2=\hbox{\box#1\unhbox#2}}\def\pusht#1\to#2{\setbox#%
2=\hbox{\unhbox#2\box#1}}

\newif\ifallowhorizmap\allowhorizmaptrue\newif\ifallowvertmap
\allowvertmapfalse\newif\ifincommdiag\incommdiagfalse

\def\diagram{\hbox\bgroup$\vcenter\bgroup\startbombverticallist
\incommdiagtrue\baselineskip\DiagramCellHeight\lineskip\z@\lineskiplimit\z@
\mathsurround\z@\tabskip\z@\let\\\diagcr\allowhorizmaptrue\allowvertmaptrue
\halign\bgroup\lcdtempl##\rcdtempl&&\lcdtempl##\rcdtempl\cr}\def\enddiagram{%
\crcr\egroup\reformatmatrix\egroup$\egroup}

\def\lcdtempl{\futurelet\thefirsttoken\dolcdtempl}\newif\ifemptycell\def
\dolcdtempl{\ifx\thefirsttoken\rcdtempl\then\hskip1sp plus 1fil \emptycelltrue
\else\hfil$\emptycellfalse\objectstyle\fi}\def\rcdtempl{\ifemptycell\else$%
\hfil\fi}\def\diagcr{\cr} \def\across#1{\span\omit\mscount=#1 \loop\ifnum
\mscount>2 \spAn\repeat\ignorespaces}\def\spAn{\relax\span\omit\advance
\mscount by -1}

\def\CellSize{\afterassignment\cdhttowd\DiagramCellHeight}\def\cdhttowd{%
\DiagramCellWidth\DiagramCellHeight}\def\MapsAbut{\MapShortFall\z@}

\newcount\cdvdl\newcount\cdvdr\newcount\cdvd\newcount\cdbfb\newcount\cdbfr
\newcount\cdbfl\newcount\cdvdr\newcount\cdvdl\newcount\cdvd

\def\reformatmatrix{\bombparameters\cdvdl=\insc@unt\cdvdr=\cdvdl\cdbfb=%
\boxc@unt\advance\cdbfb1 \cdbfr=\cdbfb\setbox1=\vbox{}\dimen2=\z@\loop\setbox
0=\lastbox\ifhbox0 \dimen1=\lastskip\unskip\dimen5=\ht0 \advance\dimen5 \dimen
1 \dimen4=\dp0 \penalty1 \reformatrow\unpenalty\ht4=\dimen5 \dp4=\dimen4 \ht3%
\z@\dp3\z@\setbox1=\vbox{\box4 \nointerlineskip\box3 \nointerlineskip\unvbox1%
}\dimen2=\dimen1 \repeat\unvbox1}

\newif\ifcontinuerow

\def\reformatrow{\cdbfl=\cdbfr\noindent\unhbox0 \loopa\unskip\setbox\cdbfl=%
\lastbox\ifhbox\cdbfl\advance\cdbfl1\repeat\par\unskip\dimen6=2%
\DiagramCellWidth\dimen7=-\DiagramCellWidth\setbox3=\hbox{}\setbox4=\hbox{}%
\setbox7=\box\voidb@x\cdvd=\cdvdl\continuerowtrue\loopa\advance\cdvd-1
\adjustcells\ifcontinuerow\advance\dimen6\wd\cdbfl\cdda=.5\dimen6 \ifdim\cdda
<\DiagramCellWidth\then\dimen6\DiagramCellWidth\advance\dimen6-\cdda
\nopendvert\cdda\DiagramCellWidth\fi\advance\dimen7\cdda\dimen6=\wd\cdbfl
\reformatcell\advance\cdbfl-1 \repeat\advance\dimen7.5\dimen6 \outHarrow} \def
\adjustcells{\ifnum\cdbfr>\cdbfl\then\ifnum\cdvdr>\cdvd\then\continuerowfalse
\else\setbox\cdbfl=\hbox to\wd\cdvd{\lcdtempl\VonH{}\rcdtempl}\fi\else\ifnum
\cdvdr>\cdvd\then\advance\cdvdr-1 \setbox\cdvd=\vbox{}\wd\cdvd=\wd\cdbfl\dp
\cdvd=\dp1 \fi\fi}

\def\reformatcell{\sethorizhtdp\noindent\unhbox\cdbfl\skip0=\lastskip\unskip
\par\ifcase\prevgraf\reformatempty\or\reformatobject\else\reformatcomplex\fi
\unskip}\def\reformatobject{\setbox6=\lastbox\unskip\vadjdon6\outVarrow
\setbox6=\hbox{\unhbox6}\advance\dimen7-.5\wd6 \outHarrow\dimen7=-.5\wd6
\pusht6\to4}\newcount\globnum

\def\reformatcomplex{\setbox6=\lastbox\unskip\setbox9=\lastbox\unskip\setbox9%
=\hbox{\unhbox9 \skip0=\lastskip\unskip\global\globnum\lastpenalty\hskip\skip
0 }\advance\globnum9999 \ifcase\globnum\reformathoriz\or\reformatpile\or
\reformatHonV\or\reformatVonH\or\reformatvert\or\reformatHmeetV\fi}

\def\reformatempty{\vpassdon\ifdim\skip0>\z@\then\hpassdon\else\ifvoid2 \then
\else\advance\dimen7-.5\dimen0 \cdda=\wd2\advance\cdda.5\dimen0\wd2=\cdda\fi
\fi}\def\VonH{\doVonH6}\def\doVonH#1{\cdna-999#1\futurelet\thenexttoken
\dooVonH}\def\dooVonH{\let\next\relax\sethorizhtdp\ifallowhorizmap
\ifallowvertmap\then\ifx\thenexttoken[\then\let\next\VonHstrut\else
\sethorizhtdp\dimen0\MapBreadth\let\next\VonHnostrut\fi\else\mapctxterr\fi
\else\mapctxterr\fi\next}\def\VonHstrut[#1]{\setbox0=\hbox{$#1$}\dimen0\wd0%
\dimen8\ht0\dimen9\dp0 \VonHnostrut}\def\VonHnostrut{\setbox0=\hbox{}\ht0=%
\dimen8\dp0=\dimen9\wd0=.5\dimen0 \copy0\penalty\cdna\box0 \allowhorizmapfalse
\allowvertmapfalse}\def\reformatHonV{\hpassdon\doreformatHonV}\def
\reformatHmeetV{\dimen@=\wd9 \advance\dimen7-\wd9 \outHarrow\setbox6=\hbox{%
\unhbox6}\dimen7-\wd6 \advance\dimen@\wd6 \setbox6=\hbox to\dimen@{\hss}%
\pusht6\to4\doreformatHonV}\def\doreformatHonV{\setbox9=\hbox{\unhbox9 \unskip
\unpenalty\global\setbox\globbox=\lastbox}\vadjdon\globbox\outVarrow}\def
\reformatVonH{\vpassdon\advance\dimen7-\wd9 \outHarrow\setbox6=\hbox{\unhbox6%
}\dimen7=-\wd6 \setbox6=\hbox{\kern\wd9 \kern\wd6}\pusht6\to4}\def\hpassdon{}%
\def\vpassdon{\dimen@=\dp\cdvd\advance\dimen@\dimen4 \advance\dimen@\dimen5
\dp\cdvd=\dimen@\nopendvert}\def\vadjdon#1{\dimen8=\ht#1 \dimen9=\dp#1 }

\def\HorizontalMap#1#2#3#4#5{\sethorizhtdp\setbox1=\makeharrowpart{#1}\def
\arrowfillera{#2}\def\arrowfillerb{#4}\setbox5=\makeharrowpart{#5}\ifx
\arrowfillera\justhorizline\then\def\arra{\hrule\horizhtdp}\def\kea{\kern-0.%
01em}\let\arrstruthtdp\horizhtdp\else\def\kea{\kern-0.15em}\setbox2=\hbox{%
\kea${\arrowfillera}$\kea}\def\arra{\copy2}\def\arrstruthtdp{height\ht2 depth%
\dp2 }\fi\ifx\arrowfillerb\justhorizline\then\def\arrb{\hrule\horizhtdp}\def
\keb{kern-0.01em}\ifx\arrowfillera\empty\then\let\arrstruthtdp\horizhtdp\fi
\else\def\keb{\kern-0.15em}\setbox4=\hbox{\keb${\arrowfillerb}$\keb}\def\arrb
{\copy4}\ifx\arrowfilera\empty\then\def\arrstruthtdp{height\ht4 depth\dp4 }%
\fi\fi\setbox3=\makeharrowpart{{#3}\vrule width\z@\arrstruthtdp}%
\ifallowhorizmap\then\let\execmap\execHorizontalMap\else\let\execmap
\mapctxterr\fi\allowhorizmapfalse\gettwoargs}\def\makeharrowpart#1{\hbox{%
\mathsurround\z@\edef\next{#1}\ifx\next\empty\else$\mkern-1.5mu{\next}\mkern-%
1.5mu$\fi}}\def\justhorizline{-}

\def\execHorizontalMap{\dimen0=\wd6 \ifdim\dimen0<\wd7\then\dimen0=\wd7\fi
\dimen3=\wd3 \ifdim\dimen0<2em\then\dimen0=2em\fi\skip2=.5\dimen0
\ifincommdiag plus 1fill\fi minus\z@\advance\skip2-.5\dimen3 \skip4=\skip2
\advance\skip2-\wd1 \advance\skip4-\wd5 \kern\MapShortFall\box1 \xleaders
\arra\hskip\skip2 \vbox{\lineskiplimit\maxdimen\lineskip.5ex \ifhbox6 \hbox to%
\dimen3 {\hss\box6\hss}\fi\vtop{\box3 \ifhbox7 \hbox to\dimen3 {\hss\box7\hss
}\fi}}\ifincommdiag\kern-.5\dimen3\penalty-9999\null\kern.5\dimen3\fi
\xleaders\arrb\hskip\skip4 \box5 \kern\MapShortFall}

\def\reformathoriz{\vadjdon6\outVarrow\ifvoid7\else\mapclasherr\fi\setbox2=%
\box9 \wd2=\dimen7 \dimen7=\z@\setbox7=\box6 }

\def\resetharrowpart#1#2{\ifvoid#1\then\ifdim#2=\z@\else\setbox4=\hbox{%
\unhbox4\kern#2}\fi\else\ifhbox#1\then\setbox#1=\hbox to#2{\unhbox#1}\else
\widenpile#1\fi\pusht#1\to4\fi}\def\outHarrow{\resetharrowpart2{\wd2}\pusht2%
\to4\resetharrowpart7{\dimen7}\pusht7\to4\dimen7=\z@}

\def\pile#1{{\incommdiagtrue\let\pile\innerpile\allowvertmapfalse
\allowhorizmaptrue\baselineskip.5\PileSpacing\lineskip\z@\lineskiplimit\z@
\mathsurround\z@\tabskip\z@\let\\\pilecr\vcenter{\halign{\hfil$##$\hfil\cr#1
\crcr}}}\ifincommdiag\then\ifallowhorizmap\then\penalty-9998
\allowvertmapfalse\allowhorizmapfalse\else\mapctxterr\fi\fi}\def\pilecr{\cr}%
\def\innerpile#1{\noalign{\halign{\hfil$##$\hfil\cr#1 \crcr}}}

\def\reformatpile{\vadjdon9\outVarrow\ifvoid7\else\mapclasherr\fi\penalty1
\setbox9=\hbox{\unhbox9 \unskip\unpenalty\setbox9=\lastbox\unhbox9 \global
\setbox\globbox=\lastbox}\unvbox\globbox\setbox9=\vbox{}\setbox7=\vbox{}%
\loopb\setbox6=\lastbox\ifhbox6 \skip3=\lastskip\unskip\splitpilerow\repeat
\unpenalty\setbox9=\hbox{$\vcenter{\unvbox9}$}\setbox2=\box9 \dimen7=\z@}\def
\pilestrut{\vrule height\dimen0 depth\dimen3 width\z@}\def\splitpilerow{%
\dimen0=\ht6 \dimen3=\dp6 \noindent\unhbox6\unskip\setbox6=\lastbox\unskip
\unhbox6\par\setbox6=\lastbox\unskip\ifcase\prevgraf\or\setbox6=\hbox\tozpt{%
\hss\unhbox6\hss}\ht6=\dimen0 \dp6=\dimen3 \setbox9=\vbox{\vskip\skip3 \hbox
to\dimen7{\hfil\box6}\nointerlineskip\unvbox9}\setbox7=\vbox{\vskip\skip3
\hbox{\pilestrut\hfil}\nointerlineskip\unvbox7}\or\setbox7=\vbox{\vskip\skip3
\hbox{\pilestrut\unhbox6}\nointerlineskip\unvbox7}\setbox6=\lastbox\unskip
\setbox9=\vbox{\vskip\skip3 \hbox to\dimen7{\pilestrut\unhbox6}%
\nointerlineskip\unvbox9}\fi\unskip}

\def\widenpile#1{\setbox#1=\hbox{$\vcenter{\unvbox#1 \setbox8=\vbox{}\loopb
\setbox9=\lastbox\ifhbox9 \skip3=\lastskip\unskip\setbox8=\vbox{\vskip\skip3
\hbox to\dimen7{\unhbox9}\nointerlineskip\unvbox8}\repeat\unvbox8 }$}}

\def\justverticalline{|}\def\makevarrowpart#1{\hbox to\MapBreadth{\hss$\kern
\MapBreadth{#1}$\hss}}\def\VerticalMap#1#2#3#4#5{\setbox1=\makevarrowpart{#1}%
\def\arrowfillera{#2}\setbox3=\makevarrowpart{#3}\def\arrowfillerb{#4}\setbox
5=\makevarrowpart{#5}\ifx\arrowfillera\justverticalline\then\def\arra{\vrule
width\MapBreadth}\def\kea{\kern-0.05ex}\else\def\kea{\kern-0.35ex}\setbox2=%
\vbox{\kea\makevarrowpart\arrowfillera\kea}\def\arra{\copy2}\fi\ifx
\arrowfillerb\justverticalline\then\def\arrb{\vrule width\MapBreadth}\def\keb
{\kern-0.05ex}\else\def\keb{\kern-0.35ex}\setbox4=\vbox{\keb\makevarrowpart
\arrowfillerb\keb}\def\arrb{\copy4}\fi\ifallowvertmap\then\let\execmap
\execVerticalMap\else\let\execmap\mapctxterr\fi\allowhorizmapfalse\gettwoargs
}

\def\execVerticalMap{\setbox3=\makevarrowpart{\box3}\setbox0=\hbox{}\ht0=\ht3
\dp0\z@\ht3\z@\box6 \setbox8=\vtop spread2ex{\offinterlineskip\box3 \xleaders
\arrb\vfill\box5 \kern\MapShortFall}\dp8=\z@\box8 \kern-\MapBreadth\setbox8=%
\vbox spread2ex{\offinterlineskip\kern\MapShortFall\box1 \xleaders\arra\vfill
\box0}\ht8=\z@\box8 \ifincommdiag\then\kern-.5\MapBreadth\penalty-9995 \null
\kern.5\MapBreadth\fi\box7\hfil}

\newcount\colno\newdimen\cdda\newbox\globbox\def\reformatvert{\setbox6=\hbox{%
\unhbox6}\cdda=\wd6 \dimen3=\dp\cdvd\advance\dimen3\dimen4 \setbox\cdvd=\hbox
{}\colno=\prevgraf\advance\colno-2 \loopb\setbox9=\hbox{\unhbox9 \unskip
\unpenalty\dimen7=\lastkern\unkern\global\setbox\globbox=\lastbox\advance
\dimen7\wd\globbox\advance\dimen7\lastkern\unkern\setbox9=\lastbox\vtop to%
\dimen3{\unvbox9}\kern\dimen7 }\ifnum\colno>0 \ifdim\wd9<\PileSpacing\then
\setbox9=\hbox to\PileSpacing{\unhbox9}\fi\dimen0=\wd9 \advance\dimen0-\wd
\globbox\setbox\cdvd=\hbox{\kern\dimen0 \box\globbox\unhbox\cdvd}\pushh9\to6%
\advance\colno-1 \setbox9=\lastbox\unskip\repeat\advance\dimen7-.5\wd6
\advance\dimen7.5\cdda\advance\dimen7-\wd9 \outHarrow\dimen7=-.5\wd6 \advance
\dimen7-.5\cdda\pusht9\to4\pusht6\to4\nopendvert\dimen@=\dimen6\advance
\dimen@-\wd\cdvd\advance\dimen@-\wd\globbox\divide\dimen@2 \setbox\cdvd=\hbox
{\kern\dimen@\box\globbox\unhbox\cdvd\kern\dimen@}\dimen8=\dp\cdvd\advance
\dimen8\dimen5 \dp\cdvd=\dimen8 \ht\cdvd=\z@}

\def\outVarrow{\ifhbox\cdvd\then\deepenbox\cdvd\pusht\cdvd\to3\else
\nopendvert\fi\dimen3=\dimen5 \advance\dimen3-\dimen8 \setbox\cdvd=\vbox{%
\vfil}\dp\cdvd=\dimen3} \def\nopendvert{\setbox3=\hbox{\unhbox3\kern\dimen6}}%
\def\deepenbox\cdvd{\setbox\cdvd=\hbox{\dimen3=\dimen4 \advance\dimen3-\dimen
9 \setbox6=\hbox{}\ht6=\dimen3 \dp6=-\dimen3 \dimen0=\dp\cdvd\advance\dimen0%
\dimen3 \unhbox\cdvd\dimen3=\lastkern\unkern\setbox8=\hbox{\kern\dimen3}%
\loopb\setbox9=\lastbox\ifvbox9 \setbox9=\vtop to\dimen0{\copy6
\nointerlineskip\unvbox9 }\dimen3=\lastkern\unkern\setbox8=\hbox{\kern\dimen3%
\box9\unhbox8}\repeat\unhbox8 }}

\newif\ifPositiveGradient\PositiveGradienttrue\newif\ifClimbing\Climbingtrue
\newcount\DiagonalChoice\DiagonalChoice1 \newcount\lineno\newcount\rowno
\newcount\charno\def\laf{\afterassignment\xlaf\charno='}\def\xlaf{\hbox{%
\tenln\char\charno}}\def\lah{\afterassignment\xlah\charno='}\def\xlah{\hbox{%
\tenln\char\charno}}\def\makedarrowpart#1{\hbox{\mathsurround\z@${#1}$}}\def
\lad{\afterassignment\xlad\charno='}\def\xlad{\setbox2=\xlaf\setbox0=\hbox to%
.5\wd2{$\hss\ldot\hss$}\ht0=.25\ht2 \dp0=\ht0 \hbox{\mv-\ht0\copy0 \mv\ht0%
\box0}}

\def\DiagonalMap#1#2#3#4#5{\ifPositiveGradient\then\let\mv\raise\else\let\mv
\lower\fi\setbox2=\makedarrowpart{#2}\setbox1=\makedarrowpart{#1}\setbox4=%
\makedarrowpart{#4}\setbox5=\makedarrowpart{#5}\setbox3=\makedarrowpart{#3}%
\let\execmap\execDiagonalLine\gettwoargs}

\def\makeline#1(#2,#3;#4){\hbox{\dimen1=#2\relax\dimen2=#3\relax\dimen5=#4%
\relax\vrule height\dimen5 depth\z@ width\z@\setbox8=\hbox to\dimen1{\tenln#1%
\hss}\cdna=\dimen5 \divide\cdna\dimen2 \ifnum\cdna=0 \then\box8 \else\dimen4=%
\dimen5 \advance\dimen4-\dimen2 \divide\dimen4\cdna\dimen3=\dimen1 \cdnb=%
\dimen2 \divide\cdnb1000 \divide\dimen3\cdnb\cdnb=\dimen4 \divide\cdnb1000
\multiply\dimen3\cdnb\dimen6\dimen1 \advance\dimen6-\dimen3 \cdnb=0
\ifPositiveGradient\then\dimen7\z@\else\dimen7\cdna\dimen4 \multiply\dimen4-1
\fi\loop\raise\dimen7\copy8 \ifnum\cdnb<\cdna\hskip-\dimen6 \advance\cdnb1
\advance\dimen7\dimen4 \repeat\fi}}\newdimen\objectheight\objectheight1.5ex

\def\execDiagonalLine{\setbox0=\hbox\tozpt{\cdna=\xcoord\cdnb=\ycoord\dimen8=%
\wd2 \dimen9=\ht2 \dimen0=\cdnb\DiagramCellHeight\advance\dimen0-2%
\MapShortFall\advance\dimen0-\objectheight\setbox2=\makeline\box2(\dimen8,%
\dimen9;.5\dimen0)\setbox4=\makeline\box4(\dimen8,\dimen9;.5\dimen0)\dimen0=2%
\wd2 \advance\dimen0-\cdna\DiagramCellWidth\advance\dimen0 2\DiagramCellWidth
\dimen2\DiagramCellHeight\advance\dimen2-\MapShortFall\dimen1\dimen2 \advance
\dimen1-\ht1 \advance\dimen2-\ht2 \dimen6=\dimen2 \advance\dimen6.25\dimen8
\dimen3\dimen2 \advance\dimen3-\ht3 \dimen4=\dimen2 \dimen7=\dimen2 \advance
\dimen4-\ht4 \advance\dimen7-\ht7 \advance\dimen7-.25\dimen8
\ifPositiveGradient\then\hss\raise\dimen4\hbox{\rlap{\box5}\box4}\llap{\raise
\dimen6\box6\kern.25\dimen9}\else\kern-.5\dimen0 \rlap{\raise\dimen1\box1}%
\raise\dimen2\box2 \llap{\raise\dimen7\box7\kern.25\dimen9}\fi\raise\dimen3%
\hbox\tozpt{\hss\box3\hss}\ifPositiveGradient\then\rlap{\kern.25\dimen9\raise
\dimen7\box7}\raise\dimen2\box2\llap{\raise\dimen1\box1}\kern-.5\dimen0 \else
\rlap{\kern.25\dimen9\raise\dimen6\box6}\raise\dimen4\hbox{\box4\llap{\box5}}%
\hss\fi}\ht0\z@\dp0\z@\box0}

\def
\NorthEast{\PositiveGradienttrue\Climbingtrue\DiagonalChoice1 }

\newif\ifmoremapargs\def\gettwoargs{\setbox7=\box\voidb@x\setbox6=\box
\voidb@x\moremapargstrue\def\whichlabel{6}\def\xcoord{2}\def\ycoord{2}\def
\contgetarg{\def\whichlabel{7}\ifmoremapargs\then\let\next\getanarg\let
\contgetarg\execmap\else\let\next\execmap\fi\next}\getanarg}\def\getanarg{%
\futurelet\thenexttoken\switcharg}\def\getlabel#1#2#3{\setbox#1=\hbox{$%
\labelstyle\>{#3}\>$}\dimen0=\ht#1\advance\dimen0 .4ex\ht#1=\dimen0 \dimen0=%
\dp#1\advance\dimen0 .4ex\dp#1=\dimen0 \contgetarg}\def\eatspacerepeat{%
\afterassignment\getanarg\let\junk= }\def\catcase#1:{{\ifcat\noexpand
\thenexttoken#1\then\global\let\xcase\docase\fi}\xcase}\def\tokcase#1:{{\ifx
\thenexttoken#1\then\global\let\xcase\docase\fi}\xcase}\def\default:{\docase}%
\def\docase#1\esac#2\esacs{#1}\def\skipcase#1\esac{}\def\getcoordsrepeat(#1,#%
2){\def\xcoord{#1}\def\ycoord{#2}\getanarg}\let\esacs\relax\def\switcharg{%
\global\let\xcase\skipcase\catcase{&}:\moremapargsfalse\contgetarg\esac
\catcase\bgroup:\getlabel\whichlabel-\esac\catcase^:\getlabel6\esac\catcase_:%
\getlabel7\esac\tokcase{~}:\getlabel3\esac\tokcase(:\getcoordsrepeat\esac
\catcase{ }:\eatspacerepeat\esac\default:\moremapargsfalse\contgetarg\esac
\esacs}

\catcode`\@=\cdna

\title{Elimination Theory for  Solvable Polynomial Algebras and Their Free Modules}

\author{Huishi Li\thanks{e-mail: huishipp@yahoo.com}\\
{\small Department of Applied Mathematics, College of Information Science and Technology}\\
{\small Hainan University,  Haikou 570228, China}}

\date{}

\begin{document} 
\maketitle
\begin{center}
\begin{minipage}{135mm}
{\small {\bf Abstract.} Let $K$ be a field, and $A=K[a_1,\ldots
,a_n]$ a solvable polynomial algebra in the sense of [K-RW, {\it J. Symbolic
Comput.}, 9(1990), 1--26].  Based on the Gr\"obner basis theory for $A$ and for free modules over $A$,  an elimination theory for  left ideals of $A$ and  an elimination theory for submodules of free $A$-modules are established.}
\end{minipage}\end{center} {\parindent=0pt\vskip 6pt

\vskip
6pt

}

\vskip .5truecm

\def\hang{\hangindent\parindent}
\def\textindent#1{\indent\llap{#1\enspace}\ignorespaces}
\def\re{\par\hang\textindent}

\def\v5{\vskip .5truecm}\def\QED{\hfill{$\Box$}}\def\hang{\hangindent\parindent}
\def\textindent#1{\indent\llap{#1\enspace}\ignorespaces}
\def\item{\par\hang\textindent}
\def \r{\rightarrow}\def\OV#1{\overline {#1}}
\def\normalbaselines{\baselineskip 24pt\lineskip 4pt\lineskiplimit 4pt}
\def\mapdown#1{\llap{$\vcenter {\hbox {$\scriptstyle #1$}}$}
                                \Bigg\downarrow}
\def\mapdownr#1{\Bigg\downarrow\rlap{$\vcenter{\hbox
                                    {$\scriptstyle #1$}}$}}
\def\mapright#1#2{\smash{\mathop{\longrightarrow}\limits^{#1}_{#2}}}
\def\NZ{\mathbb{N}}\def\mapleft#1#2{\smash{\mathop{\longleftarrow}\limits^{#1}_{#2}}}

\def\LH{{\bf LH}}\def\LM{{\bf LM}}\def\LT{{\bf
LT}}\def\KX{K\langle X\rangle} \def\KS{K\langle X\rangle}
\def\B{{\cal B}} \def\LC{{\bf LC}} \def\G{{\cal G}} \def\FRAC#1#2{\displaystyle{\frac{#1}{#2}}}
\def\SUM^#1_#2{\displaystyle{\sum^{#1}_{#2}}} \def\O{{\cal O}}  \def\J{{\bf J}}\def\BE{\B (e)}
\def\PRCVE{\prec_{\varepsilon\hbox{-}gr}}\def\BV{\B (\varepsilon )}\def\PRCEGR{\prec_{e\hbox{\rm -}gr}}

\def\KS{K\langle X\rangle}
\def\LR{\langle X\rangle}\def\T#1{\widetilde #1}
\def\HL{{\rm LH}}\def\HY{\hbox{\hskip .03truecm -\hskip .03truecm}}\def\F{{\cal F}}

\def\hang{\hangindent\parindent}
\def\textindent#1{\indent\llap{#1\enspace}\ignorespaces}
\def\re{\par\hang\textindent}

\def\hang{\hangindent\parindent}
\def\textindent#1{\indent\llap{#1\enspace}\ignorespaces}
\def\item{\par\hang\textindent}
\def\normalbaselines{\baselineskip 24pt\lineskip 4pt\lineskiplimit 4pt}

\def\v5{\vskip .5truecm}
\def\QED{\hfill{$\Box$}}

\def\mapdown#1{\llap{$\vcenter {\hbox {$\scriptstyle #1$}}$}
                                \Bigg\downarrow}
\def\mapdownr#1{\Bigg\downarrow\rlap{$\vcenter{\hbox
                                    {$\scriptstyle #1$}}$}}

\def\mapright#1#2{\smash{\mathop{\longrightarrow}\limits^{#1}_{#2}}}
\def\mapleft#1#2{\smash{\mathop{\longleftarrow}\limits^{#1}_{#2}}}
\def \r{\rightarrow}

\def\FRAC#1#2{\displaystyle{\frac{#1}{#2}}}
\def\SUM^#1_#2{\displaystyle{\sum^{#1}_{#2}}}

\def\OV#1{\overline {#1}}
\def\T#1{\widetilde {#1}}

\def\MB{\mathbb{B}}

\def\LH{{\bf LH}}
\def\LM{{\bf LM}}
\def\LT{{\bf LT}}
\def\LC{{\bf LC}}

\def\B{{\cal B}}
\def\O{{\cal O}}
\def\G{{\cal G}}
\def\F{{\cal F}}
\def\S{{\cal S}}

\def\VF{\varphi}
\def\BE{\B (e)}
\def\VE{\varepsilon}
\def\PRCVE{\prec_{\varepsilon\textrm{\tiny -}gr}}
\def\BV{\B (\varepsilon )}

\def\PRCEGR{\prec_{e\textrm{\tiny -}gr}}
\def\PRECSVE{\prec_{s\textrm{\tiny -}\varepsilon}}
\def\PRECEVE{\prec_{e\textrm{\tiny -}\VE}}

\def\KX{K\langle X\rangle}
\def\GR{Gr\"obner~}


This note is a complement to [Li5] so that the \GR basis theory for solvable polynomial algebras and free modules over such algebras is getting more complete.

\section*{1. Solvable Polynomial algebras}

{\bf 1.1. Definition} ([K-RW], [LW]) Let $A=K[a_1,\ldots ,a_n]$ be a
finitely generated $K$-algebra. Suppose that  $A$ has the PBW
$K$-basis $\B=\{a^{\alpha}=a_1^{\alpha_1}\cdots
a_n^{\alpha_n}~|~\alpha =(\alpha_1,\ldots ,\alpha_n)\in\NZ^n\}$,
and that $\prec$ is a (two-sided) monomial ordering on $\B$. If for
all $a^{\alpha}=a_1^{\alpha_1}\cdots a_n^{\alpha_n}$,
$a^{\beta}=a_1^{\beta_1}\cdots a^{\beta_n}_n\in\B$, the following
holds:
$$\begin{array}{rcl} a^{\alpha}a^{\beta}&=&\lambda_{\alpha ,\beta}a^{\alpha +\beta}+f_{\alpha ,\beta},\\
&{~}&\hbox{where}~\lambda_{\alpha ,\beta}\in K^*,~a^{\alpha
+\beta}=a_1^{\alpha_1+\beta_1}\cdots a_n^{\alpha_n+\beta_n}~\hbox{and}~f_{\alpha ,\beta}\in K\hbox{-span}\B\\
&{~}&\hbox{such that either}~f_{\alpha ,\beta}=0~\hbox{or}~\LM
(f_{\alpha ,\beta})\prec a^{\alpha +\beta},\end{array}\leqno{(\hbox{S})}$$ where $\LM
(f_{\alpha ,\beta})$ stands for the leading monomial of $f_{\alpha
,\beta}$ with respect to $\prec$, then $A$ is called a {\it solvable
polynomial algebra}.\par

Usually $(\B ,\prec )$ is referred to an {\it admissible system} of $A$.\v5

The next proposition provides us with a constructive characterization of solvable polynomial algebras.{\parindent=0pt\v5

{\bf 1.2. Proposition}  [Li3, Theorem 2.1] Let $A=K[a_1,\ldots
,a_n]$ be a finitely generated $K$-algebra, and let $\KS =K\langle
X_1,\ldots ,X_n\rangle$ be the free $K$-algebras with the standard
$K$-basis $\mathbb{B}=\{ 1\}\cup\{X_{i_1}\cdots X_{i_s}~|~X_{i_j}\in
X,~s\ge 1\}$. The following two statements are equivalent:}\par

(i) $A$ is a solvable polynomial algebra in the sense of Definition
2.1.\par

(ii) $A\cong \OV A=\KS /I$ via the $K$-algebra epimorphism $\pi_1$:
$\KS \r A$ with $\pi_1(X_i)=a_i$, $1\le i\le n$, $I=$ Ker$\pi_1$,
satisfying  {\parindent=1.85truecm

\item{(a)} with respect to some monomial ordering $\prec_{_X}$ on $\mathbb{B}$, the ideal $I$ has a
finite Gr\"obner basis $G$ and the reduced Gr\"obner basis of $I$ is
of the form
$$\G =\left\{ g_{ji}=X_jX_i-\lambda_{ji}X_iX_j-F_{ji}
~\left |~\begin{array}{l} \LM (g_{ji})=X_jX_i,\\ 1\le i<j\le
n\end{array}\right. \right\} $$ where $\lambda_{ji}\in K^*$,
$\mu^{ji}_q\in K$, and
$F_{ji}=\sum_q\mu^{ji}_qX_1^{\alpha_{1q}}X_2^{\alpha_{2q}}\cdots
X_n^{\alpha_{nq}}$ with $(\alpha_{1q},\alpha_{2q},\ldots
,\alpha_{nq})\in\NZ^n$, thereby $\B =\{ \OV X_1^{\alpha_1}\OV
X_2^{\alpha_2}\cdots \OV X_n^{\alpha_n}~|~$ $\alpha_j\in\NZ\}$ forms
a PBW $K$-basis for $\OV A$, where each $\OV X_i$ denotes the coset
of $I$ represented by $X_i$ in $\OV A$; and

\item{(b)} there is a (two-sided) monomial ordering
$\prec$ on $\B$ such that $\LM (\OV{F}_{ji})\prec \OV X_i\OV X_j$
whenever $\OV F_{ji}\ne 0$, where $\OV F_{ji}=\sum_q\mu^{ji}_q\OV
X_1^{\alpha_{1i}}\OV X_2^{\alpha_{2i}}\cdots \OV X_n^{\alpha_{ni}}$,
$1\le i<j\le n$. \par}\QED\v5

We refer to [K-RW] and [Li1] for the  propositions 1.3 -- 1.5 stated below.{\parindent=0pt\v5

{\bf 1.3. Proposition}   Let $A=K[a_1,\ldots ,a_n]$ be a solvable
polynomial algebra with   admissible system $(\B ,\prec )$. The
following statements hold.\par

(i) If $f,g\in A$ with $\LM (f)=a^{\alpha}$, $\LM (g)=a^{\beta}$,
then
$$\begin{array}{rcl} \LM (fg)&=&\LM (\LM (f)\LM (g))\\
&=&\LM (a^{\alpha}a^{\beta})\\
&=&a^{\alpha +\beta}\\
&=&\LM (a^{\beta}a^{\alpha})\\
&=&\LM (\LM (g)\LM (f))\\
&=&\LM (gf).\end{array}$$\par

(ii) $A$ is a domain, that is, $A$ has no (left and right) divisors
of zero.\QED\v5

{\bf 1.4. Proposition}  Let $A=K[a_1,\ldots ,a_n]$ be a solvable
polynomial algebra with   admissible system $(\B ,\prec )$. Consider the polynomial ring $A[t]$ over  $A$ where $t$ is a commuting variable $t$, i.e., $tf=ft$ for all $f\in A$. Then $A[t]$ is a free $A$-module with the $A$-basis $\{ t^q~|~q\in\NZ\}$ and it is also a $K$-algebra with the generating set $\{ a_1,\ldots ,a_n,t\}$ and the  PBW basis
$$\B (t)=\{ a^{\alpha}t^q~|~a^{\alpha}\in\B ,~q\in\NZ\}.$$
Furthermore,  $\B (t)$ can be equipped with a monomial ordering $\prec_t$ subject to the rule: for $a^{\alpha}t^{\ell}$, $a^{\beta}t^q\in\B (t)$,
$$a^{\alpha}t^{\ell}\prec_ta^{\beta}t^q\Leftrightarrow\left\{\begin{array}{l}
\ell <q\\
\hbox{or}\\
\ell =q~\hbox{and}~ a^{\alpha}\prec a^{\beta},\end{array}\right.$$
which turns $A[t]$ into a solvable polynomial algebra.\QED\v5

{\bf 1.5. Proposition}  Let $A_1=K[a_1,\ldots ,a_n]$ and
$A_2=K[b_1,\ldots ,b_m]$ be solvable polynomial algebras with
admissible systems $(\B_1,\prec_1)$ and $(\B_2,\prec_2)$
respectively, where $\B_1=\{ a^{\alpha}=a_1^{\alpha_1}\cdots a_n^{\alpha_n}~|~\alpha =(\alpha_1,\ldots ,\alpha_n)\in\NZ^n\}$ and $\B_2=\{ b^{\beta}=b_1^{\alpha_1}\cdots b_m^{\alpha_m}~|~\beta =(\beta_1,\ldots ,\beta_m)\in\NZ^m\}$  Then the tensor product $A=A_1\otimes_KA_2$ is a solvable polynomial
algebra with the admissible system $(\B ,\prec)$, where
$\B=\{a^{\alpha}\otimes
b^{\beta}~|~a^{\alpha}\in\B_1,~b^{\beta}\in\B_2\}$, while $\prec$ is
defined on $\B$ subject to the rule: for $a^{\alpha}\otimes
b^{\beta}$, $a^{\gamma}\otimes b^{\eta}\in\B$,
$$a^{\alpha}\otimes
b^{\beta}\prec a^{\gamma}\otimes
b^{\eta}\Leftrightarrow\left\{\begin{array}{l}
a^{\alpha}\prec_1a^{\gamma}\\
\hbox{or}\\
a^{\alpha}=a^{\gamma}~\hbox{and}~b^{\beta}\prec_2b^{\eta}.\end{array}\right.$$}\v5

In [K-RW] the \GR basis theory for (one-sided, two-sided) ideals of solvable polynomial algebras were developed. In particular,  every left ideal $N$ of a solvable polynomial algebra $A$ has a finite \GR basis $\G$ in the sense that if $f\in N$ and $f\ne 0$, then there is a $g\in\G$ such that $\LM (g)|_{_{\rm L}}\LM (f)$, i.e., there is a monomial  $a^{\gamma}\in\B$ such that $\LM (f)=\LM (a^{\gamma}\LM (g))$, thereby $A$ is a Noetherian domain because the same is true for every right ideal;  moreover, there is a noncommutative version of Buchberger's algorithm that computes a \GR basis for every finitely generated left ideal of $N$.

\section*{2. Elimination Orderings and Elimination in Left Ideals}

By introducing the elimination ordering on a solvable polynomial algebra $A$ with respect to a subspace  instead of a subalgebra of $A$ (comparing with the elimination theory based on the \GR basis theory for commutative polynomial ideals [AL2, Section 2.3]),  showing that every solvable polynomial algebra has an elimination ordering, and applying an elimination lemma given in [Li4] to solvable polynomial algebras, in this section we establish an elimination theory for left ideals of solvable polynomial algebras. A more general elimination theory for submodules of free modules over solvable polynomial algebras will be developed in  Section 4.\v5

 We start with the following{\parindent=0pt\v5

{\bf 2.1. Definition}  Let $A=K[a_1,\ldots ,a_n]$ be a solvable polynomial algebra with admissible system $(\B ,\prec)$. Given a nonempty subset $S\subset\B$,
let $V(S)=K$-span$S$. If the monomial ordering $\prec$ on $\B$ is such that
$$f\in A~\hbox{and}~\LM (f)\preceq a^{\alpha}~\hbox{for some}~a^{\alpha}\in S~\hbox{implies}~f\in V(S),$$
then it is referred to as an {\it elimination ordering} with respect to $V(S)$.\v5

{\bf Remark} The idea of defining an elimination monomial ordering as described in Definition 2.1 above has been enlightened by [Li1, Section 4 -- Section 7 of CH.V] and proposed in [Li4] for binomial skew polynomial rings and general solvable polynomial polynomial algebras, but not yet explicitly written down as a formal  definition in loc. cit.\v5

{\bf Example} (1) Let $A=K[a_1,\ldots ,a_n]$ be a finitely generated $K$-algebra such that
for all $1\le i<j\le n$,
$$\begin{array}{rcl} a_ja_i&=&\lambda_{ji}a_ia_j+\sum\lambda_{\alpha}x_{j_1}^{\alpha_1}x_{j_2}^{\alpha_2}\cdots x_{j_t}^{\alpha_t}+\mu_{ji},\\
&{~}&\hbox{where}~j_1<j_2<\cdots <j_t<i,~t\le n-1,\\
&{~}&\quad\quad\quad\lambda_{ji}\in K^*,~\lambda_{\alpha},\mu_{ji}\in K.\end{array}$$
Then $A$ is a solvable polynomial algebra with  the  monomial ordering $\prec_{lex}$ which is the
lexicographic ordering such that
$$a_n\prec_{lex}a_{n-1}\prec_{lex}\cdots\prec_{lex}a_1.$$
 Moreover, $\prec_{lex}$ is an elimination ordering with respect to $V(S)=K$-span$S$,  where  for $2\le r\le n$, $$S=\{a_r^{\alpha_r}a_{r+1}^{\alpha_{r+1}}\cdots a_n^{\alpha_n}~|~(\alpha_r,\ldots ,\alpha_n)\in\NZ^{n-r+1}\} ,$$
This example may be applied to  the Weyl algebra $A_n(K)$, and more generally, to the additive analogue of the Weyl algebra $A=A_n(q_1,\ldots ,q_n)=K[x_1,\ldots ,x_n,y_1,\ldots ,y_n]$  in the sense of [Kur] and [JBS], many Ore extensions, many skew polynomial algebras, many operator algebras listed in [K-RW] and [Li1], and enveloping algebras $U(g)$ of many finite dimensional Lie algebras (such as Heisenberg algebra).}\v5

From the literature (e.g. [K-RW], [Li1]) one may also find out  (or try to establish) more elimination orderings over  other solvable polynomial algebras. Below let us present two more examples which may also be  of independent interest.{\parindent=0pt\v5

{\bf Example} (2) Let $A=K[a_1,\ldots ,a_n]$ be a solvable
polynomial algebra with   admissible system $(\B ,\prec )$. Consider the polynomial extension $A[t]$ of $A$ by a commuting variable $t$, i.e., $tf=ft$ for all $f\in A$. Then, as a $K$-algebra generated by $\{ a_1,\ldots ,a_n,t\}$, $A[t]$ is a solvable polynomial algebra with the admissible system $(\B (t) ,\prec_t)$, where $\B (t)=\{ a^{\alpha}t^q~|~a^{\alpha}\in\B ,~q\in\NZ\}$ is the PBW basis of $A[t]$ and the monomial ordering $\prec_t$ is the one defined in Proposition 1.4 such that $a^{\alpha}\prec_t a^{\beta}t^q$ for all $a^{\alpha},a^{\beta}\in\B$ and $t^q\ne 1$.  So, one sees that $\prec_t$ is an elimination ordering  with respect to $V(S)$, where $S=\B$ which is the PBW basis of $A$. \v5

{\bf Example} (3) Let $A_1=K[a_1,\ldots ,a_n]$ and
$A_2=K[b_1,\ldots ,b_m]$ be solvable polynomial algebras with
admissible systems $(\B_1,\prec_1)$ and $(\B_2,\prec_2)$
respectively. Then, as a $K$-algebra generated by $\{a_i\otimes 1,~1\otimes b_j~|~1\le i\le n,~1\le j\le m\}$, the tensor product $A=A_1\otimes_KA_2$ is a solvable polynomial
algebra with the PBW basis $\B=\{a^{\alpha}\otimes
b^{\beta}~|~a^{\alpha}\in\B_1,~b^{\beta}\in\B_2\}$ and the monomial ordering $\prec$  defined in Proposition 1.5 such that $1\otimes b^{\beta}\prec a^{\alpha}\otimes b^{\gamma}$ for all $b^{\beta},b^{\gamma}\in\B_2$ and $a^{\alpha}\in\B_1-\{ 1\}$. So one sees that $\prec$ is  an elimination ordering with respect to $V(S)$, where $S=\{ 1\otimes b^{\beta}~|~b^{\beta}\in\B_2\}\subset\B$.}\v5

Indeed, for a solvable polynomial algebra $A=K[a_1,\ldots ,a_n]$, the proposition below tells us that with respect to any proper subset $U$ of $\{ a_1,\ldots ,a_n\}$,  $A$ has an elimination monomial ordering in the sense of Definition 2.1.{\parindent=0pt\v5

{\bf 2.2. Proposition}   Let $A=K[a_1,\ldots ,a_n]$ be a solvable polynomial algebra with   admissible system $(\B ,\prec )$. For any subset $U=\{ a_{i_1},a_{i_2}\ldots ,a_{i_m}\}\subset\{a_1,a_2,\ldots ,a_n\}$ with $i_1<i_2<\cdots <i_m$ and $m<n$, if we consider the subset
$$S=\left\{\left. a^{\beta}= a_{i_1}^{\beta_1}a_{i_2}^{\beta_2}\cdots a_{i_m}^{\beta_{m}}~\right |~\beta =(\beta_1,\beta_2,\ldots ,\beta_{m})\in\NZ^{m}\right\}$$ of $\B$, then the given monomial ordering $\prec$ on $\B$ gives rise to an elimination ordering $\lessdot$ on $\B$ with respect to $V(S)=K$-span$S$. Moreover the restriction of $\lessdot$ on $S$ coincides with the restriction of $\prec$ on $S$ (that is $\prec$). \vskip 6pt

{\bf Proof} For the given subset $U$, putting $U^c=\{a_1,\ldots ,a_n\}-U$,  we may write $U^c=\{a_{j_1},$ $a_{j_2},$ $\ldots ,$ $a_{j_{n-m}}\}$ such that $j_1<j_2<\cdots <j_{n-m}$, and thus we may put $$S^c=\left\{\left. a^{\alpha}=a_{j_1}^{\alpha_1}a_{j_2}^{\alpha_2}\cdots a_{j_{n-m}}^{\alpha_{n-m}}~\right |~\alpha =(\alpha_1,\alpha_2,\ldots ,\alpha_{n-m})\in\NZ^{n-m}\right\} .$$
Note that by the definition of a solvable polynomial algebra (Definition 1.1), any $a^{\gamma}\in\B$ may be uniquely written as $a^{\gamma}=\LM (a^{\alpha}a^{\beta})$ for some $a^{\alpha}\in S^c$ and $a^{\beta}\in S$. Note also that if $a^{\alpha (1)},a^{\alpha}\in S^c$, $a^{\beta (1)},a^{\beta}\in S$, then it follows from Proposition 1.3 that
$$\begin{array}{l} \LM (a^{\alpha (1)}a^{\alpha})=a^{\alpha (1)+\alpha}\in S^c,\\
\LM (a^{\beta (1)}a^{\beta})=a^{\beta (1)+\beta}\in S.\end{array}$$
Thereby if $a^{\eta}\in\B$ and $a^{\eta}=\LM (a^{\alpha (1)}a^{\beta (1)})$, then with $a^{\gamma}=\LM (a^{\alpha}a^{\beta})$, again by Proposition 1.3 we have
$$\begin{array}{rcl} \LM (a^{\eta}a^{\gamma})&=&\LM \left(\LM (a^{\alpha (1)}a^{\beta (1)})\LM (a^{\alpha}a^{\beta})\right )\\
&=&a^{(\alpha (1)+\alpha )+(\beta (1)+\beta )}\\
&=&\LM \left (a^{\alpha (1)+\alpha }a^{\beta (1)+\beta }\right )\\
&=&\LM \left (\LM (a^{\alpha (1)}a^{\alpha})\LM (a^{\beta (1)}a^{\beta})\right ).\end{array}$$
So, if we define on $\B$ a new ordering $\lessdot$ subject to the rule: for $a^{\gamma (1)}=\LM (a^{\alpha (1)}a^{\beta (1)})$, $a^{\gamma (2)}=\LM (a^{\alpha (2)}a^{\beta (2)})\in\B$,
$$a^{\gamma (1)}\lessdot a^{\gamma (2)}\Leftrightarrow\left\{\begin{array}{l} a^{\alpha (1)}\prec a^{\alpha (2)}\\
\hbox{or}\\
a^{\alpha (1)}= a^{\alpha (2)}~\hbox{and}~a^{\beta (1)}\prec a^{\beta (2)},\end{array}\right.$$
then by referring to Definition 1.1, Proposition 1.3 and the formulas derived above, one checks that $\lessdot$ is a monomial ordering on $\B$ such that
$$a^{\beta}\lessdot a^{\alpha}, \quad a^{\beta}\in S,~a^{\alpha}\in S^c~\hbox{and}~a^{\alpha}\ne 1.$$
Furthermore, by the definition of $\lessdot$ and Definition 2.1, we see that $\lessdot$ is an elimination ordering on $\B$ with respect to $V(S)=K$-span$S$. Finally, the assertion concerning  the restriction of $\lessdot$ on $S$ follows from the construction of $\lessdot$.\QED}\v5

With an elimination ordering in the sense of Definition 2.1, the elimination principle via \GR bases of  left ideals in a solvable polynomial algebra is embodied by the following{\parindent=0pt \v5

{\bf 2.3. Theorem}  Let $A=K[a_1,\ldots a_s]$ be a solvable polynomial algebra with admissible system $(\B,\prec )$, where $\B$ is the PBW basis of $A$ and for a certain subset $S\subset\B$, $\prec$ is an elimination monomial ordering on $\B$ with respect to  $V(S)=K$-span$S$, and let $N$ be a left ideal of $A$.  If $\G$ is a \GR basis of $N$ with respect to the elimination ordering $\prec$ on $\B$, then the following statements hold.}\par

(i) If $f\in N\cap V(S)$ and $f\ne 0$, then there is a $g\in\G\cap V(S)\subset N\cap V(S)$ such that $\LM (g)~|_{_{\rm L}}\LM (f)$ and, if this is the case, then there is an $a^{\gamma}\in S$ such that $\LM (f)=\LM (a^{\gamma}\LM (g))$ and thus $$f=\lambda^{-1}\mu a^{\gamma}g+f_1~\hbox{with}~\LM (f)=\LM (a^{\gamma}g),~\LM (f_1)\prec\LM (f),\leqno{(*)}$$  where $\lambda =\LC (g)$ and $\mu =\LC (f)$. \par

(ii) Let $D\subset A$ be a subalgebra of $A$ (conventionally $D$ and $A$ have the same identity element 1), and suppose that $S$ is a $K$-basis for $D$. Then with respect to the restriction of $\prec$ on $S$, every nonzero $f\in N\cap D$ has an expression
$$\begin{array}{rcl}f&=&\sum_{i,j}\nu_{i,j}a^{\gamma (i)}g_j~\hbox{with}~\nu_{i,j}\in K^*,~a^{\gamma (i)}\in S\subset D,\\
&{~}&g_j\in\G\cap V(S)=\G\cap D\subset N\cap D,\\
&{~}&\hbox{such that}~\LM (a^{\gamma (i)}g_j)\preceq\LM (f)~\hbox{ for all appearing}~(i,j).\end{array}$${\parindent=0pt\vskip 6pt

{\bf Proof} (i) Let $f\in N\cap V(S)$ be a nonzero element. As $\G$ is a \GR basis of $N$ in $A$  with respect to $\prec$, it follows  that there is a $g\in\G$  such that $\LM (g)|_{_{\rm L}}\LM (f)$ in $A$, thereby $\LM (g)\preceq\LM (f)$. Note also  $f\in V(S)$ and thus $\LM (f)=a^{\alpha}$ for some $a^{\alpha}\in S$. It turns out that $\LM (g)\preceq a^{\alpha}$. Since  $\prec$ is an elimination monomial ordering   on $\B$ with respect to $V(S)$, it follows from Definition 2.1 that  $g\in \G\cap V(S)\subset N\cap V(S)$. Turning back to $\LM (g)|_{_{\rm L}}\LM (f)$, there is an $a^{\gamma}\in\B$ such that $\LM (f)=\LM (a^{\gamma}\LM (g))$. But by Proposition 1.3(i) we thus have $\LM (f)=\LM (\LM (g)a^{\gamma})$, namely $a^{\gamma}|_{_{\rm L}}\LM (f)$. So, similarly from $a^{\gamma}\preceq\LM (f)$ we derive $a^{\gamma}\in S$. Hence, the desired expression $(*)$ is obtained.}  \par

(ii) By the assumption on $S$, we have $V(S)=D$. So, it follows from  (i) that $N\cap D\ne \{ 0\}$ implies $\G\cap D\ne\emptyset$. Thus for a nonzero element $f\in N\cap D$,  the expression $(*)$ in (i) entails $f_1=f -\lambda^{-1}\mu a^{\gamma}g\in N\cap D$ for some $g\in \G\cap D$. If $f_1\ne 0$, then repeat this division procedure on $f_1$ by $\G\cap D$ and so on. As $\prec$ is a well-ordering,  after a finite number of repeating the division procedure by $\G\cap D$ we then reach the desired expression for $f$. \QED{\parindent=0pt\v5

{\bf Remark} We observe that in the case of Theorem 1.4.3 (ii), no matter the subalgebra $D$ of $A$ is a solvable polynomial algebra or not (with respect to the restriction of $\prec$ on $S$), $\G\cap D$ is indeed a \GR basis for the left ideal $N\cap D$ of $D$ (in the sense that every nonzero element has a \GR representation with respect to the restriction of $\prec$ on $S$).}\v5

Comparing with the elimination theorem for commutative polynomial ideals (e.g. [AL2, Theorem 2.3.4]), we now derive a similar result for left ideals of solvable polynomial algebras.
To see this, let $A=K[a_1,\ldots ,a_n]$ be a solvable polynomial algebra with   admissible system $(\B ,\prec )$ in which $\prec$ is an elimination monomial ordering on $\B$ with respect to $V(S)$ (in the sense of Definition 2.1), where for a given subset $U=\{ a_{i_1},a_{i_2}\ldots ,a_{i_m}\}\subset\{a_1,a_2,\ldots ,a_n\}$ with $i_1<i_2<\cdots <i_m$ and $m<n$, $$S=\left\{\left. a^{\beta}= a_{i_1}^{\beta_1}a_{i_2}^{\beta_2}\cdots a_{i_m}^{\beta_{m}}~\right |~\beta=(\beta_1,\beta_2,\ldots ,\beta_{m})\in\NZ^{m}\right\}\subset\B $$  and $V(S)=K$-span$S$ (note that by Proposition 2.2, such an elimination ordering $\prec$ exists).{\parindent=0pt\v5

{\bf 2.4. Theorem}     With the notation above, let $N$ be a left ideal of $A$, and let $\G$ be a \GR basis of $N$ with respect to the elimination ordering $\prec$ on $A$.  Consider  the subalgebra $K[U]$ of $A$ generated by $U$, and suppose that the subset  $S$ of $\B$ forms a $K$-basis for $K[U]$. If $N\cap K[U]\ne\{ 0\}$,  then $\G\cap K[U]$ is a Gr\"obner basis  of the left ideal $N\cap K[U]$ in $K[U]$ with respect to $\prec_{_{S}}$, where $\prec_{_{S}}$ is the restriction of $\prec$ on $S$.{\parindent=0pt\vskip 6pt

{\bf Proof} This follows immediately from Theorem 2.3.\par\QED\v5

{\bf Remark} With the notation used above,  let us emphasise that in Theorem 2.4 the assumption  made on $S$ is necessary, though $V(S)\subset K[U]$. This is  because elements of $\G$ are linear combinations of elements in $\B$, while in principle, getting elements of $\G\cap K[U]$ amounts to  eliminating generators out of $U$ from certain elements $g$ of $\G$ such that $g \in \G\cap V(S)$. So, if $f\in N\cap K[U]-V(S)$, then the given elimination ordering cannot help to reach the result of Theorem 2.3(i). But in general $S$ may not necessarily a $K$-basis for $K[U]$. For instance, one may look at the $q$-Heisenberg algebra ${\bf h}_n(q)$ whcih is
 generated over a field $K$ by the set of elements $\{
x_i,y_i,z_i~|~i=1,...,n\}$ subject to the relations:
$$\begin{array}{ll}
x_ix_j=x_jx_i,~y_iy_j=y_jy_i,~z_jz_i=z_iz_j,& 1\le i<j\le n,\\
x_iz_i=qz_ix_i,&1\le i\le n,\\
z_iy_i=qy_iz_i,&1\le i\le n,\\
x_iy_i=q^{-1}y_ix_i+z_i,&1\le i\le n,\\
x_iy_j=y_jx_i,~x_iz_j=z_jx_i,~y_iz_j=z_jy_i,&i\ne j,\\
\end{array}$$
where $q\in K^*$, and consider the subalgebra $K[U]$ generated by the subset $U=\{ x_i,y_i~|~1\le i\le n\}$. Or else let us look at the solvable polynomial algebra $A=K[a_1,a_2,a_3]$ given in [Li3],  where the generators of $A$ satisfies $a_2a_1=a_1a_2$, $a_3a_1=\lambda a_1a_3+\mu a_2^2a_3+f(a_2)$ with $f(a_2)\in K$-span$\{ 1,a_2,a_2^2,\ldots ,a_2^6\}$, $a_3a_2=a_2a_3$,  and $A$ has the PBW basis  $\mathcal{B}' =\{ a^{\alpha}=a_1^{\alpha_1}a_2^{\alpha_2}a_3^{\alpha_3}~|~\alpha
=(\alpha_1,\alpha_2,\alpha_3)\in\NZ^3\}$ and the monomial ordering $\prec_{lex}$ on $\B '$ such that $a_3\prec_{lex}a_2\prec_{lex}a_1$. If $K[U]$ is the subalgebra of $A$ generated by the subset $U=\{ a_1,a_3\}\subset\{ a_1,a_2,a_3\}$, and $S=\{a_1^{\alpha_1}a_3^{\alpha_3}~|~(\alpha_1,\alpha_3)\in\NZ^2\}\subset\B '$, then  $a_3^{\alpha_3}a_1^{\alpha_1}=a_1^{\alpha_1}a_3^{\alpha_3}+\sum\mu_{\gamma_2 ,\eta_3}a_2^{\gamma_2}a_3^{\eta_3}\not\in V(S)=K$-span$S$ entails that $S$ is not a $K$-basis of $K[U]$.}}\v5

As to the problem mentioned in the remark above, actually we have the following easily verified fact.{\parindent=0pt

\v5
{\bf 2.5. Proposition}   Let $A=K[a_1,\ldots ,a_n]$ be a solvable polynomial algebra with   admissible system $(\B ,\prec )$ and let $K[U]$ be the subalgebra of $A$ generated by a subset $U=\{ a_{i_1},a_{i_2}\ldots ,a_{i_m}\}\subset\{a_1,a_2,\ldots ,a_n\}$ with $i_1<i_2<\cdots <i_m$ and $m<n$.
Then the subset
$$S=\left\{\left. a^{\beta}= a_{i_1}^{\beta_1}a_{i_2}^{\beta_2}\cdots a_{i_m}^{\beta_{m}}~\right |~\beta =(\beta_1,\beta_2,\ldots ,\beta_{m})\in\NZ^{m}\right\}$$ is a $K$-basis for $K[U]$ if and only if $K[U]$ is a solvable polynomial algebra with the PBW basis $S$ and the admissible system $(S,\prec_S)$, where $\prec_S$ is the restriction of $\prec$ on $S$.\QED}\v5

In contrast to Theorem 2.4, we will see in the next section that without any assumption on $S$, an elimination theorem  does hold true for submodules of a free module $L=\oplus_{i=1}^sAe_i$ (with $s\ge 2$) over a solvable polynomial algebra $A$.\v5

As an application of the above Example (2) and Theorem 2.3, we next derive  a result for determining the generating set of intersection of two left ideals via elimination in left ideals. To see this, Let $A=K[a_1,\ldots ,a_n]$ be a solvable polynomial algebra with   admissible system $(\B ,\prec )$, and let $N_1$, $N_2$ be two left ideals of $A$. Considering the polynomial extension $A[t]$ of $A$ by a commuting variable $t$, let
$$N=\sum_{u\in N_1}A[t]tu+\sum_{v\in N_2}A[t](1-t)v$$
be the left ideal of $A[t]$ generated by $\{ tu, (1-t)v~|~u\in N_1,~v\in N_2\}$. {\parindent=0pt\v5

{\bf 2.6. Proposition}  With the notation fixed above, the following statements hold.}\par

(i) $N\cap A=N_1\cap N_2$.\par

(ii) If $\G$ is a \GR basis of $N$ with respect to the monomial ordering $\prec_t$ on $A[t]$ (as constructed in Proposition 1.4), then $\G\cap A$ is a \GR basis for the left ideal $N\cap A=N_1\cap N_2$ of $A$ with respect to the monomial ordering $\prec$ on $A$. {\parindent=0pt\vskip 6pt

{\bf Proof} (i) If $f\in N_1\cap N_2$,  then $f=tf+(1-t)f\in N$, showing $N_1\cap N_2\subseteq N\cap A$. On the other hand, if $f\in A\cap N$, say $f=tg+(1-t)h$ with $g\in\sum_{u\in N_1}A[t]u$ and $h\in\sum_{v\in N_2}A[t]v$, then, considering the algebra homomorphism
$$\begin{array}{cccc} \varphi_1:&A[t]&\mapright{}{}&A\\
&t&\mapsto&0\end{array}$$
from which we see that $f\in A$ entails $f=\varphi_1(f)=\varphi_1(h)\in N_2$, and considering the algebra homomorphism
$$\begin{array}{cccc} \varphi_2:&A[t]&\mapright{}{}&A\\
&t&\mapsto&1\end{array}$$
from which we see that $f\in A$ entails $f=\varphi_2(f)=\varphi_2(g)\in N_1$. Hence $f\in N_1\cap N_2$ and consequently $A\cap N\subseteq N_1\cap N_2$. Combining inclusions of both directions we conclude  $N_1\cap N_2=A\cap N$.}\par

(ii)   By Example (2) above, the monomial ordering $\prec_t$ on $A[t]$ is an elimination ordering  with respect to $V(S)$, where $S=\B\subset\B (t)$.  Noticing that $A$ is a subalgebra of $A[t]$ with the PBW basis $\B$ and the restriction of $\prec_t$ on $\B$ coincides with the monomial ordering $\prec$ on $\B$, it follows from Theorem 2.3(ii) that if $\G$ is a \GR basis of $N$ with respect to the monomial ordering $\prec_t$ on $A[t]$,  then $\G\cap A$ is a \GR basis for the left ideal $N\cap A=N_1\cap N_2$ of $A$ with respect to $\prec$.\par\QED\v5

Let $A=K[a_1,\ldots ,a_n]$ be a solvable polynomial algebra with the PBW basis $\B$. Then it follows from Proposition 2.2 that for any given subset $U=\{ a_{i_1},a_{i_2}\ldots ,a_{i_m}\}\subset\{a_1,a_2,\ldots ,a_n\}$ with $i_1<i_2<\cdots <i_m$ and $m<n$,  $A$ has an elimination monomial ordering $\prec$ on $\B$ with respect to $V(S)$ (in the sense of Definition 2.1), where $S=\left\{\left. a^{\beta}= a_{i_1}^{\beta_1}a_{i_2}^{\beta_2}\cdots a_{i_m}^{\beta_{m}}~\right |~\beta=(\beta_1,\beta_2,\ldots ,\beta_{m})\in\NZ^{m}\right\}\subset\B $  and $V(S)=K$-span$S$.
We observe that  even if $S$ is not the PBW basis of the subalgebra $K[U]$ generated by $U$, or equivalently, even if $K[U]$ is not a solvable polynomial algebra with the admissible system $(S,\prec)$,  if for a certain true reason concerning a left ideal $N$ of $A$ we need to take some elements from $N\cap V(S)$, then fortunately Theorem 2.3(i) tells us that the elimination ordering $\prec$ may still enable us to take out {\it desired elements} from $G\cap V(S)$ {\it once we know} $N\cap V(S)\ne\{0 \}$, where  $\G$ is a \GR basis of the left ideal $N$ with respect to $\prec$. With this merit of Proposition 2.2 and Theorem 2.3, our purpose below is to determine and realize the elimination property for left ideals of $A$ (in the sense of [Li4, Lemma 3.1]) by means of \GR bases.\v5

Recall from [Li1, Section 2 of CH.III] that a solvable polynomial  algebra $A=K[a_1,\ldots ,a_n]$ is called a  {\it quadric solvable polynomial algebra} if for all $1\le i<j\le n$,
 $$a_ja_i=\lambda_{ij}a_ia_j+\sum_{q\le \ell}\mu_{q\ell}a_qa_{\ell}+\sum_kc_ka_k+c,\quad \lambda_{ij}\in K^*,\mu_{q\ell},c_k,c\in K,$$
and $A$ admits a graded monomial ordering $\prec_{gr}$
with $d(a_i)=1$ for $1\le i\le n$. For instance, Weyl algebras and enveloping algebras of finite dimensional Lie algebras are typical examples of such algebras. Since  $A$  has the PBW basis $\B =\{ a^{\alpha}= a_{1}^{\alpha_1}\cdots
a_{n}^{\alpha_n}~|~\alpha =(\alpha_1,\ldots ,\alpha_n)\in\NZ^n\}$, a subset
$U=\{a_{i_1},\ldots ,a_{i_r}\}\subset\{a_1,\ldots ,a_n\}$ with
$i_1<i_2<\cdots <i_r$ is said to be {\it weakly independent
modulo a left ideal  $N$ of $A$} if $N\cap V(S)=\{0\}$, where
$$S=\left\{\left. a^{\alpha}= a_{i_1}^{\alpha_1}\cdots
a_{i_r}^{\alpha_r}~\right |~\alpha =(\alpha_1,...,\alpha_r)
\in\NZ^r\right\}\subset\B$$
and $V(S)=K$-span$S$
With the weak independence of $U$ (mod $N$) and  a double
filtered-graded transfer trick,  the strategy of computing dim${\cal
V}(I)$ proposed by [KW] was adapted in [Li1, CH.V] to
compute the Gelfand-Kirillov dimension GK.dim$(A/N)$, and
consequently the following results were established: \vskip 6pt
{\parindent=.44truecm

\item{$\bullet$}  [Li1, CH.V, Theorem  7.4]  Let $A$ be a quadric solvable polynomial algebra, and  $N$ a nonzero left ideal of $A$. Then
$$\begin{array}{rcl}\hbox{GK.dim}(A/N)
&=&\hbox{degree of the Hilbert polynomial of}~A/N\\
&=&\max\left\{|U|~\left |~\begin{array}{l} U\subset\{a_1,\ldots ,a_n\}~\hbox{is}\\
\hbox{weakly independent}~(\hbox{mod}~N)\end{array}\right.\right\},\end{array}$$  which can be
algorithmically computed  via a Gr\"obner basis of $N$; moreover,
$$\hbox{GK.dim}(A/N)< n=\hbox{GK.dim}A;$$

\item{$\bullet$} [Li1, CH.V, Lemma 7.5]  If GK.dim$A/N=d$, then $V(S)\cap N\ne \{ 0\}$
for every subset $U=\{ a_{i_1},...,a_{i_{d+1}}\}\subset \{
a_1,...,a_n\}$ with $i_1<i_2<\cdots <i_{d+1}$, where $S=\{ a^{\alpha}=a_{i_1}^{\alpha_1}\cdots
a_{i_{d+1}}^{\alpha_{d+1}}~|~\alpha =(\alpha_1,...,\alpha_{d+1})
\in\NZ^{d+1}\}\subset\B$ and $V(S)=K$-span$S$.\vskip
6pt}

Note that the class of quadric solvable polynomial algebras studied
in [Li1, CH.III, CH.V] covers not only Weyl algebras and enveloping algebras of Lie algebras, but also more  Ore extensions, skew polynomial algebras, and operator algebras. Enlightened by the
automatic proving of multivariate identities over operator algebras
([PWZ], [Ch], [CS]),   more general  $\partial$-finiteness and
$\partial$-holonomicity for modules over quadric solvable polynomial
algebras have been introduced and preliminarily studied  in [Li1. CH.VII]
where [Li1, CH.V, Lemma 7.5] has played a key role.\v5

Turning to an arbitrary solvable polynomial algebra  in the sense of Definition 1.1,
we have the following{\parindent=0pt\v5

{\bf 2.7. Theorem}  Let $A=K[a_1,\ldots ,a_n]$ be a solvable
polynomial algebra with admissible system $(\B ,\prec )$. For every $1\le d< n-1$ and the subset $U=\{ a_{i_1},...,a_{i_{d+1}}\}\subset$ $\{ a_1,...,a_n\}$ with
$i_1<i_2<\cdots <i_{d+1}$, let $S=\{ a^{\alpha}=a_{i_1}^{\alpha_1}\cdots
a_{i_{d+1}}^{\alpha_{d+1}}~|~\alpha =(\alpha_1,...,\alpha_{d+1})
\in\NZ^{d+1}\}\subset\B$ and $V(S)=K$-span$S$. Then the following statements hold.}\par

(i) GK.dim$A=n$ ($=$ the number of generators of
$A$), and GK.dim$(A/N)<n$ holds for every nonzero left ideal $N$ of $A$.\par

(ii) If GK.dim$(A/N)=d$, then  $V(S)\cap N\ne
\{ 0\}$, in particular for every $U=\{
a_1,\ldots a_s\}$ with $d+1\le s\le n-1$ we have $V(S)\cap
N\ne \{ 0\}$.  In orther words, Elimination lemma  holds true (in the sense of [Li4, Definition 3.2]) for every nonzero left ideal $N$ of $A$.{\parindent=0pt\vskip 6pt

{\bf Proof} We prove both the assertions (i) and (ii) together. First note that every solvable polynomial algebra $A$
has the PBW basis $\B$ by Definition 1.1. Moreover, it follows from
Proposition 1.2 and [Li2, Example 1 of Section 5.3] that  GK.dim$A=n$ ($=$ the number of generators of
$A$). As also we know that $A$ is a domain by Proposition 1.3(ii). Hence [Li4, Lemma
3.3]   entails that GK.dim$(A/N)<$ GK.dim$A=n$ holds for every
nonzero left ideal of $A$. Therefore, we conclude that Elimination
lemma  holds true (in the sense of [Li4, Definition 3.2]) for every
nonzero left ideal $N$ of $A$.  \QED}\v5

Furthermore, from [K-RW] we know that a noncommutative Buchberger algorithm  works effectively for every solvable polynomial algebra $A$, that is, if a finite
generating set of a left ideal $N$ of $A$ is given (note that $A$ is
Noetherian), then running the noncommutative Buchberger algorithm
with respect to a monomial ordering $\prec$ will produce a finite
Gr\"obner basis $\G$ for $N$.  Thus it follows from Theorem 2.7 that we have the next{\parindent=0pt\v5

{\bf 2.8. Theorem} {old 1.4.8) Let $A=K[a_1,\ldots, a_n]$ be a solvable polynomial algebra with admissible system $(\B ,\prec )$, and let $N$ be a nonzero left ideal of $A$ with GK.dim$(A/N)=d$. For any  given subset $U=\{ a_{i_1},a_{i_2},\ldots ,a_{i_{d+1}}\}\subset\{a_1,\ldots ,a_n\}$ with $i_1<i_2,\ldots <i_{d+1}$ and the subset $S=\{ a^{\alpha}=a_{i_1}^{\alpha_1}a_{i_2}^{\alpha_2}\cdots a_{d+1}^{\alpha_{d+1}}~|~\alpha =(\alpha_1,\alpha_2,\ldots ,\alpha_{d+1})\in\NZ^{d+1}\}\subset\B ,$
let $\prec$ be an elimination ordering on $\B$ with respect to $V(S)=K$-span$S$ (in the sense of Definition 1.4.1). If $\G$ is a Gr\"obner basis of $N$ with respect to $\prec$,  then $\G\cap V(S)\ne\emptyset$.\vskip 6pt

{\bf Proof} By Theorem 2.7, we have $N\cap V(S)\ne \{ 0\}$. It follows from  Theorem 2.3(i) that  $\G\cap V(S)\ne\emptyset$.\QED\v5

{\bf 2.9. Corollary}  Let $A=K[a_1,\ldots, a_n]$ be a solvable polynomial algebra with admissible system $(\B ,\prec )$, and let $N$ be a nonzero left ideal of $A$ with GK.dim$(A/N)=d$. Then for any  given subset $U=\{ a_{i_1},a_{i_2},\ldots ,a_{i_{d+1}}\}\subset\{a_1,\ldots ,a_n\}$ with $i_1<i_2,\ldots <i_{d+1}$ and the subset $S=\{ a^{\alpha}=a_{i_1}^{\alpha_1}a_{i_2}^{\alpha_2}\cdots a_{d+1}^{\alpha_{d+1}}~|~\alpha =(\alpha_1,\alpha_2,\ldots ,\alpha_{d+1})\in\NZ^{d+1}\}\subset\B ,$  getting a nonzero element from the intersection $N\cap V(S)$ may be realized  in (at most) two  steps:}}\par

{\bf Step 1.} If the given $\prec$ is not an elimination monomial ordering on $\B$ with respect to $V(S)$, then,   with the subset $S$ in use, construct the elimination ordering $\lessdot$ with respect to $V(S)$ as described in Proposition 2.2;\par

{\bf Step 2.} Run the noncommutative Buchberger algorithm for solvable polynomial algebras to produce a Gr\"obner basis $\G$ for $N$ with respect to $\lessdot$, and then take out the desired nonzero element from $\G\cap V(U)$.\v5

\section*{3. \GR Bases of Submodules in Free Modules}

Let $A=K[a_1,\ldots ,a_n]$ be a solvable polynomial algebra with
admissible system $(\B ,\prec )$ in the sense of Definition 1.1.3,
where $\B =\{ a^{\alpha}=a_1^{\alpha_1}\cdots
a_n^{\alpha_n}~|~\alpha =(\alpha_1,\ldots ,\alpha_n)\in\NZ^n\}$ is
the PBW basis of $A$ and $\prec$ is a monomial ordering on $\B$,
and let $L=\oplus_{i=1}^sAe_i$ be a free {\it left} $A$-module with the
$A$-basis $\{ e_1,\ldots ,e_s\}$. Then $L$ has the $K$-basis
$$\BE =\{ a^{\alpha}e_i~|~a^{\alpha}\in\B ,~1\le i\le s\} .$$ For
convenience, elements of $\BE$ are also referred to as {\it
monomials} in $L$.\v5

If $\prec_{e}$ is a total ordering on $\BE$, and if $\xi
=\sum_{j=1}^m\lambda_ja^{\alpha (j)}e_{i_j}\in L$, where
$\lambda_j\in K^*$ and $\alpha (j)=(\alpha_{j_1},\ldots
,\alpha_{j_n})\in\NZ^n$, such that $$a^{\alpha (1)}e_{i_1}\prec_{e}
a^{\alpha (2)}e_{i_2}\prec_{e}\cdots\prec_{e} a^{\alpha
(m)}e_{i_m},$$ then by $\LM (\xi )$ we denote the {\it leading
monomial} $a^{\alpha (m)}e_{i_m}$ of $\xi $, by $\LC (\xi )$ we
denote the {\it leading coefficient} $\lambda_m$ of $\xi $,  and by
$\LT (\xi )$ we denote the {\it leading term} $\lambda_ma^{\alpha
(m)}e_{i_m}$ of $\xi$.{\parindent=0pt\v5

{\bf 3.1. Definition}   With respect to the given monomial ordering
$\prec$ on $\B$, a total ordering $\prec_{e}$ on $\BE$ is called a
{\it monomial ordering} if the following two conditions are
satisfied:}\par

(1) $a^{\alpha}e_i\prec_{e} a^{\beta}e_j$ implies  $\LM
(a^{\gamma}a^{\alpha}e_i)\prec_{e} \LM (a^{\gamma}a^{\beta}e_j)$ for
all $a^{\alpha}e_i$, $a^{\beta}e_j\in\BE$, $a^{\gamma}\in\B$;\par

(2) $a^{\alpha}\prec a^{\beta}$ implies $a^{\alpha}e_i\prec_{e}
a^{\beta}e_i$ for all $a^{\alpha},a^{\beta}\in\B$ and $1\le i\le
s$. \par

If $\prec_e$ is a monomial ordering on $\BE$, then we also say
that $\prec_e$ is a monomial ordering on the free module $L$,
and the data $(\BE ,\prec_e)$ is referred to as an {\it
admissible system} of $L$.\v5

By referring to Proposition 1.3, we record two easy but useful
facts on a monomial ordering  $\prec_{e}$ on $\BE$, as follows.{\parindent=0pt\v5

{\bf 3.2. Lemma}  (i) Every monomial ordering  $\prec_{e}$ on
$\BE$ is a well-ordering, i.e., every nonempty subset of $\BE$ has a
minimal element.} \par

(ii) If  $f\in A$ with $\LM (f)=a^{\gamma}$ and $\xi\in L$ with $\LM
(\xi )=a^{\alpha}e_i$, then
$$\LM (f\xi )=\LM (\LM (f)\LM (\xi ))=\LM (a^{\gamma}a^{\alpha}e_i)=a^{\gamma +\alpha}e_i.$$
\par\QED\v5

Actually as in the commutative case ([AL], [Eis], [KR]), any
monomial ordering $\prec$ on $\B$ may induce  two  monomial orderings on $\BE$:
$$\hbox{{\bf TOP} monomial ordering}\quad a^{\alpha}e_i\prec_{top} a^{\beta}e_j\Leftrightarrow\left\{\begin{array}{l}
a^{\alpha}\prec a^{\beta}\\
\hbox{or}\\
a^{\alpha}=a^{\beta}~
\hbox{and}~i<j;\end{array}\right.$$
$$\hbox{{\bf POT} monomial ordering}\quad  a^{\alpha}e_i\prec_{pot}
a^{\beta}e_j\Leftrightarrow\left\{\begin{array}{l} i<j\\
\hbox{or}\\
i=j~\hbox{and}~a^{\alpha}\prec a^{\beta},\end{array}\right.$$
where TOP abbreviates the phrase "term over position", while POT abbreviates the phrase "position over term".\par

Let $d(~)$ be a positive-degree function on $A$  such that $d(a_i)=m_i>0$, $1\le i\le n$, and let
$(b_1,\ldots b_s)\in\NZ^n$ be any fixed $s$-tuple. Then, by
assigning $e_j$ the degree $b_j$, $1\le j\le s$, every monomial
$a^{\alpha}e_j$ in the $K$-basis $\BE$ of $L$ is endowed with the
degree $d(a^{\alpha})+b_j$. If a
monomial ordering $\prec_{e}$ on $\BE$ satisfies
$$a^{\alpha}e_i\prec_ea^{\beta}e_j~\hbox{implies}~d(a^{\alpha})+b_i\le d(a^{\beta})+b_j,$$
then we call it a {\it graded monomial ordering} on
$\BE$ (or equivalently, a graded monomial ordering on $L$).\v5

Let  $\prec_{e}$ be a  monomial ordering on $L$, and let
$L_1=\oplus_{i=1}^mA\varepsilon_i$ be another free $A$-module with
the $A$-basis $\{ \varepsilon_1,\ldots ,\varepsilon_m\}$. Then, as
in the commutative case ([AL], [Eis], [KR]), for any given finite
subset $G=\{ g_1,\ldots ,g_m\}\subset L$,  an ordering on the
$K$-basis $\B (\varepsilon )=\{
a^{\alpha}\varepsilon_i~|~a^{\alpha}\in\B ,~1\le i\le m\}$ of $L_1$
can be defined subject to the rule: for
$a^{\alpha}\varepsilon_i,a^{\beta}\varepsilon_j\in\B (\varepsilon
)$,
$$a^{\alpha}\varepsilon_i\prec_{s\textrm{\tiny -}\varepsilon} a^{\beta}\varepsilon_j
\Leftrightarrow\left\{\begin{array}{l}
\LM (a^{\alpha}g_i)\prec_{e}\LM (a^{\beta}g_j)\\
\hbox{or}\\
\LM (a^{\alpha}g_i)=\LM (a^{\beta}g_j)~\hbox{and}~i<j.
\end{array}\right.$$
It is an exercise to check that this ordering is a monomial ordering on  $\B (\varepsilon )$.  $\prec_{s\textrm{\tiny -}\varepsilon}$ is usually referred to as the {\it Schreyer ordering} induced by $G$
with respect to $\prec_{e}$.  \v5

With respect to a given monomial ordering $\prec_{e}$ on $\BE$,
every nonzero submodule $N$ of $L$ has a finite left Gr\"obner basis
$\G =\{ g_1,\ldots ,g_m\}\subset N$ in the sense that
{\parindent=1truecm\par

\item{$\bullet$} if $\xi\in N$ and $\xi\ne 0$,
then $\LM (g_i)|_{_{\rm L}}\LM (\xi )$ for some $g_i\in\G$, i.e., there is a
monomial $a^{\gamma}\in\B$ such that $\LM (\xi )=\LM (a^{\gamma}\LM
(g_i))$, or equivalently, $\xi$ has a {\it left Gr\"obner
representation} $\xi =\sum_{i,j}\lambda_{ij}a^{\alpha (i_j)}g_j$,
where $\lambda_{ij}\in K^*$, $a^{\alpha (i_j)}\in\B$ with $\alpha
(i_j)=(\alpha_{i_{j1}},\ldots ,\alpha_{i_{jn}})\in\NZ^n$,
$g_j\in\G$,  satisfying $\LM (a^{\alpha (i_j)}g_j)\preceq_{e}\LM
(\xi)$.\par}{\parindent=0pt
Moreover, starting with any finite generating set of $N$, such a
left Gr\"obner basis $\G$ can be computed by running a
noncommutative version of the Buchberger algorithm for submodules over
solvable polynomial algebras. For more details on the \GR basis theory for modules over solvable polynomial algebras, one is referred to [K-RW], [Li1], [Lev], and [Li5]}\v5

\section*{4.  Elimination Orderings and Elimination in Submodules of Free Modules}

This section  is devoted to an elimination theory for submodules of free {\it left} modules over solvable polynomial algebras. Comparing with Section 2, we start by introducing the notion of an elimination ordering on free modules.
{\parindent=0pt\v5

{\bf 4.1. Definition}  (Compare with Definition 2.1.) Let $A=K[a_1,\ldots ,a_n]$ be a solvable polynomial algebra with admissible system $(\B ,\prec)$, and let $L=\oplus_{i=1}^sAe_i$ be a free $A$-module with admissible system $(\BE ,\prec_e)$. For a nonempty subset $S$ of $\BE =\{ a^{\alpha}e_i~|~a^{\alpha}\in\B ,~1\le i\le s\}$,  let $V(S)=K$-span$S$. If the monomial ordering  $\prec_e$ on $\BE$ is such that
$$\xi\in L~\hbox{and}~\LM (\xi )\preceq_ea^{\alpha}e_{\ell}~\hbox{for some}~a^{\alpha}e_{\ell}\in S~\hbox{implies}~\xi\in V(S),$$
then it is referred to as an {\it elimination ordering} with respect to $V(S)$. }\v5

Before working with Definition 4.1, for the convenience of examining whether the ordering  respectively defined in Example (1) -- Example (3) below is an elimination monomial ordering with respect to the given  $V(S)$, let us record an easy but useful fact concerning a special case where $S$ is the $K$-basis of a free submodule of $L$.  {\parindent=0pt\v5

{\bf 4.2. Lemma}  Let $L=\oplus_{i=1}^sAe_i$ be a free $A$-module with admissible system $(\BE ,\prec_e)$.  For any $1\le d\le s-1$ and a subset $U=\{ e_{i_1},e_{i_2},\ldots ,e_{i_d}\}\subset\{ e_1,\ldots ,e_s\}$ with $i_1<i_2<\cdots <i_d$, let  $L_U=\oplus_{j=1}^{d}Ae_{i_j}$ be the submodule of $L$ generated by $U$ and $S=\{ a^{\alpha}e_{i_j}~|~a^{\alpha}\in\B ,~e_{i_j}\in U\}$ (which is the $K$-basis of $L_U$). Then it follows from  Lemma 3.2 that
$$\LM (x^{\gamma}a^{\alpha}e_{i_j})=a^{\gamma +\alpha}e_{i_j}\in S$$
holds for all $a^{\gamma}\in\B$, $a^{\alpha}e_{i_j}\in S$.\par\QED}\v5

In the following examples, $A=K[a_1,\ldots ,a_n]$ is a solvable polynomial algebra with admissible system $(\B ,\prec )$.{\parindent=0pt\v5

{\bf Example} (1) Let $\prec_{pot}$ be the POT monomial ordering on the free $A$-module $L=\oplus_{i=1}^sAe_i$ induced by a given  monomial ordering $\prec$ on $A$ (see Section 3). Then $e_1\prec_{pot}e_2\prec_{pot}\cdots\prec_{pot}e_s$. For any $1\le d\le s-1$, let $U=\{ e_1,\ldots ,e_d\}\subset\{ e_1,\ldots ,e_n\}$, $S=\{a^{\alpha}e_{\ell}~|~a^{\alpha}\in\B ,~e_{\ell}\in U\}\subset\B (e)$  and  $V(S)=K$-span$S$ ($=L_U=\oplus_{\ell =1}^dAe_{\ell}$).  If $\xi =\lambda a^{\beta}e_t+\sum_{q,j }\lambda_{q,j }a^{\beta (q)}e_{j}\in L$ and $\LM (\xi )=a^{\beta}e_t\prec_e a^{\alpha}e_{\ell}$ for some $a^{\alpha}e_{\ell} \in S$, then $ t\le\ell\le d$ and $j\le t\le\ell\le d$, implying that $\xi\in V(S)$. This shows that the POT monomial ordering $\prec_{pot}$ on $L$ is an elimination ordering with respect to $V(S)$. Moreover, bearing in mind  Lemma 4.2, it is straightforward to check that the restriction of $\prec_{pot}$ on $S$ is a monomial ordering on $L_U$.\v5

{\bf Example} (2) Consider the free $A$-module $L=(\oplus_{i=1}^sAe_i)\oplus(\oplus_{j=1}^mA\VE_j)$. Let $\BE$ and $\BV$ denote the $K$-bases of the free submodules $L_1=\oplus_{i=1}^sAe_i$ and $L_2=\oplus_{j=1}^mA\VE_j$ respectively, then $L$ has the $K$-basis $B(e,\VE ) =\BE\cup\BV .$ If $\prec_{_{L}}$ is a  monomial ordering on $\B (e,\VE )$, then by Lemma 4.2, it is straightforward to see that the restriction of $\prec_{_{L}}$ on $\BE$, respectively on $\BV$, is a monomial ordering on $\BE$, respectively a monomial ordering on $\BV$. Now, let us define a new ordering $\lessdot$ on the $K$-basis $B(e,\VE )$ of $L$ subject to the rule: for $x,y\in \B (e,\VE )$,
$$x\lessdot y\Leftrightarrow\left\{\begin{array}{l}
x=a^{\alpha (1)}e_{i_1},y=a^{\alpha (2)}e_{i_2}\in\BE ~\hbox{and}~ a^{\alpha (1)}e_{i_1}\prec_{_L}a^{\alpha (2)}e_{i_2},\\
\hbox{or}\\
x=a^{\beta (1)}\VE_{j_1},y=a^{\beta (2)}\VE_{j_2}\in\BV ~\hbox{and}~ a^{\beta (1)}\VE_{j_1}\prec_{_L}a^{\beta (2)}\VE_{j_2},\\
\hbox{or}\\
x=a^{\beta}\VE_j\in\BV~\hbox{and}~y=a^{\alpha}e_i\in\BE .\end{array}\right.$$
Then, by using Lemma 4.2, a direct verification shows that $\lessdot$ is a monomial ordering on $\B (e,\VE )$ such that $a^{\beta}\VE_j\lessdot a^{\alpha}e_i$ for all $a^{\beta}\VE_j\in\BV$ and $a^{\alpha}e_i\in\BE$. It follows that if $\xi\in L$ and $\LM (\xi )\lessdot a^{\beta}\VE_j$ for some $a^{\beta}\VE_j\in\BV$, then $\xi\in L_2$. This shows that with   $S=\BV\subset\B (e,\VE )$ and $V(S)=K$-span$S$ ($=L_2$),  $\lessdot$ is an elimination ordering on $\B (e,\VE )$ with respect to $V(S)$. Moreover, it is clear that the restriction of $\lessdot$ on $L_1$ and $L_2$ respectively coincides with the restriction of $\prec_{_L}$ on $L_1$ and $L_2$ respectively, thereby the restriction is a monomial ordering on $L_1$ and $L_2$ respectively.}\v5

Immediately, example (2) given above motivates the next{\parindent=0pt\v5

{\bf Example} (3) Let $L_1=\oplus_{i=1}^sAe_i$ and $L_2=\oplus_{j=1}^mA\VE_j$ be free $A$-modules with admissible systems $(\B (e),\prec_{_{L_1}})$, $(\B (\VE ),\prec_{_{L_2}})$ respectively. Put $L=L_1\oplus L_2$. Then $L$ is a free $A$-module with the $K$-basis $\B (e,\VE )=\B (e)\cup\B (\VE )$.  If we define the ordering $\lessdot$ on the $K$-basis $B(e,\VE )$ of $L$ subject to the rule: for $x,y\in \B (e,\VE )$,
$$x\lessdot y\Leftrightarrow\left\{\begin{array}{l}
x=a^{\alpha (1)}e_{i_1},y=a^{\alpha (2)}e_{i_2}\in\BE ~\hbox{and}~ a^{\alpha (1)}e_{i_1}\prec_{_{L_1}}a^{\alpha (2)}e_{i_2},\\
\hbox{or}\\
x=a^{\beta (1)}\VE_{j_1},y=a^{\beta (2)}\VE_{j_2}\in\BV ~\hbox{and}~ a^{\beta (1)}\VE_{j_1}\prec_{_{L_2}}a^{\beta (2)}\VE_{j_2},\\
\hbox{or}\\
x=a^{\beta}\VE_j\in\BV~\hbox{and}~y=a^{\alpha}e_i\in\BE ,\end{array}\right.$$
then, by using Lemma 4.2,  it is straightforward to check that $\lessdot$ is a monomial ordering on $\B (e,\VE )$ , such that $a^{\beta}\VE_j\lessdot a^{\alpha}e_i$ for all $a^{\beta}\VE_j\in\BV$ and $a^{\alpha}e_i\in\BE$. It follows that if $\xi\in L$ and $\LM (\xi )\lessdot a^{\beta}\VE_j$ for some $a^{\beta}\VE_j\in\BV$, then $\xi\in L_2$. This shows that with   $S=\BV\subset\B (e,\VE )$ and $V(S)=K$-span$S$ ($=L_2$), $\lessdot$ is an elimination ordering on $\B (e,\VE )$ with respect to $V(S)$. Moreover,  the restriction of $\lessdot$ on $\BE$ and $\BV$ respectively coincides with $\prec_{_{L_1}}$ and $\prec_{_{L_2}}$ respectively, thereby the restriction is a monomial ordering on $\BE$ and $\BV$ respectively.}\v5

The next example shows that elimination orderings on free $A$-modules (in the sense of Definition 4.1) are not always like those constructed in Example (1) -- Example (3) by taking $S$ being the $K$-basis of a free $A$-submodule $L_U$ of a given free $A$-module $L$ (i.e., $L_U$ and $L$ are modules over the same ground algebra $A$).{\parindent=0pt\v5

{\bf Example} (4) Let $L=\oplus_{i=1}^sAe_i$ be a free $A$-module with admissible system $(\BE ,\prec_e)$. Considering the polynomial extension $A[t]$ of $A$ by a commuting variable $t$, then $A[t]$ is a solvable polynomial algebra with the PBW basis $\B (t)=\{ a^{\alpha}t^q~|~a^{\alpha}\in\B ,~q\in\NZ\}$ and the monomial ordering $\prec_t$ as constructed in Proposition 1.4. Let $L_t=\oplus_{i=1}^sA[t]e_i$ be the free $A[t]$-module. Then $L_t$ has the $K$-basis
$$\B (t,e)=\{ t^qa^{\alpha}e_i~|~a^{\alpha}\in\B,~q\in\NZ ,~1\le i\le s\} .$$
Note that $A$ is a subalgebra of $A[t]$, thereby $\BE\subset\B (e,t)$ and $L$ is now an $A$-submodule of $L_t$ instead of an $A[t]$-submodule of $L_t$. If we define the ordering $\prec_{t\textrm{\tiny -}e}$ on $\B (t,e)$ subject to the rule: for $t^qa^{\alpha}e_i$, $t^{\ell}a^{\beta}e_j\in\B (t,e)$,
$$t^qa^{\alpha}e_i\prec_{t\textrm{\tiny -}e} t^{\ell}a^{\beta}e_j\Leftrightarrow\left\{\begin{array}{l}
q<\ell\\
\hbox{or}\\
q=\ell~\hbox{and}~a^{\alpha}e_i\prec_e a^{\beta}e_j,\end{array}\right.$$
then it is straightforward to verify that $\prec_{t\textrm{\tiny -}e}$ is a monomial ordering on $\B (t,e)$ such that $a^{\alpha}e_i\prec_{t\textrm{\tiny -}e}t^qa^{\beta}e_j$ for all $a^{\alpha}e_i, a^{\beta}e_j\in \B (e)$ and $t^q\ne 1$. It follows that if $\xi\in L_t$ and  $\LM (\xi )\prec_{t\textrm{\tiny -}e} a^{\alpha}e_i$ for some $a^{\alpha}e_i\in\B (e)$, then $\xi\in L$. This shows that with  $S=\B (e)\subset\B (e,t)$ and $V(S)=K$-span$S$ ($=L)$, $\prec_{t\textrm{\tiny -}e}$ is an elimination ordering on $\B (t, e )$ with respect to $V(S)$. Moreover, the restriction of $\prec_{t\textrm{\tiny -}e}$ on $\B (e)$ coincides with $\prec_e$.}\v5

With an elimination ordering in the sense of Definition 4.1, we first have an analogue of Theorem 2.3 which embodies the elimination principle via \GR bases of  submodules  in a free module. To see this, Let $A=K[a_1,\ldots a_s]$ be a solvable polynomial algebra with admissible system $(\B,\prec )$, and let $L=\oplus_{i=1}^sAe_i$ be a free $A$-module with admissible system $(\BE ,\prec_e)$, where  for a certain subset $S\subset\BE$, $\prec_e$ is an elimination monomial ordering on $\BE$ with respect to  $V(S)=K$-span$S$. Furthermore, we put
$$\begin{array}{l} S_{\B}=\left\{\left. a^{\gamma}\in\B~\right |~a^{\gamma}e_i\in S~\hbox{for some}~e_i\in\{e_1,\ldots ,e_s\}\right\} ,\\
S_E=\left\{ e_j\in\{a_1,\ldots ,e_s\}~\left |~a^{\alpha}e_j\in S~\hbox{for some}~a^{\alpha}\in \B \right.\right \} .\end{array}$${\parindent=0pt \v5

{\bf 4.3. Theorem}   With the notation fixed above, considering a submodule $N$ of $L$, let $\G$ be  a \GR basis of $N$ with respect to the elimination ordering $\prec_e$ on $\BE$. Then the following statements hold.}\par

(i) If $\xi\in N\cap V(S)$ and $\xi\ne 0$, then there is a $g\in\G\cap V(S)\subset N\cap V(S)$ such that $\LM (g)~|_{_{\rm L}}\LM (\xi )$ and, if this is the case, then there is an $a^{\gamma}\in S_{\B}$ such that $\LM (\xi )=\LM (a^{\gamma}\LM (g))$ and thus $$\xi =\lambda^{-1}\mu a^{\gamma}g+\xi_1~\hbox{with}~\LM (\xi )=\LM (a^{\gamma}g),~\LM (\xi_1)\prec\LM (\xi ),\leqno{(*)}$$  where $\lambda =\LC (g)$ and $\mu =\LC (\xi )$. \par

(ii) Let $D$ be a subalgebra of $A$ (conventionally $D$ and $A$ have the same identity element 1), and let $L_D=\oplus_{e_j\in S_E}De_j$ be the free $D$-module (which is a $D$-submodule of $L$). Suppose that $S_{\B}$ is a $K$-basis for $D$. Then,
$$\left\{\begin{array}{l}L_D~\hbox{has the}~K\hbox{-basis}~\left\{\left. a^{\alpha}e_j~\right |~a^{\alpha}\in S_{\B}, ~e_j\in S_E\right\} =S\subset\B (e) ,\\
\hbox{such that}~V(S)= L_D,\end{array}\right.\leqno{(**)}$$
and with respect to the restriction $\prec_e$ on $S$, every nonzero $\xi\in N\cap L_D$ has an expression
$$\begin{array}{rcl}\xi &=&\sum_{i,j}\nu_{i,j}a^{\gamma (i)}g_j~\hbox{with}~\nu_{i,j}\in K^*,~a^{\gamma (i)}\in S_{\B}\subset D,\\
&{~}&g_j\in\G\cap V(S)=\G\cap L_D\subset N\cap L_D,\\
&{~}&\hbox{such that}~\LM (a^{\gamma (i)}g_j)\preceq_e\LM (\xi )~\hbox{ for all appearing}~(i,j).\end{array}$${\parindent=0pt\vskip 6pt

{\bf Proof} (i) Let $\xi\in N\cap V(S)$ be a nonzero element. As $\G$ is a \GR basis of $N$ in $L$  with respect to $\prec_e$, it follows   that there is a $g\in\G$  such that $\LM (g)|_{_{\rm L}}\LM (\xi )$ in $L$, thereby $\LM (g)\preceq_e\LM (\xi )$. Note also  $\xi\in V(S)$ and thus $\LM (\xi)=a^{\beta}e_i$ for some $a^{\beta}e_i\in S$. It turns out that $\LM (g)\preceq_e a^{\beta}e_i$. Since  $\prec_e$ is an elimination monomial ordering   on $\BE$ with respect to $V(S)$, it follows from Definition 4.1 that  $g\in \G\cap V(S)\subset N\cap V(S)$. Turning back to $\LM (g)|_{_{\rm L}}\LM (\xi )$, there is an $a^{\gamma}\in\B$ such that $\LM (\xi )=\LM (a^{\gamma}\LM (g))$, and if $\LM (\xi )=a^{\beta}e_i$ then the leading monomial of $g$ must be of the form $\LM (g)=a^{\alpha}e_i$. It follows from Lemma 3.2(ii) that
$$\begin{array}{rcl} a^{\beta}e_i=\LM (\xi )&=&\LM (a^{\gamma}\LM (g))\\
&=&\LM (a^{\gamma}a^{\alpha}e_i)\\
&=&a^{\gamma +\alpha}e_i\\
&=&a^{\alpha +\gamma}e_i\\
&=&\LM (a^{\alpha}a^{\gamma}e_i),\end{array}$$ showing $a^{\gamma}e_i|_{_{\rm L}}\LM (\xi )$. So, as with $g$ above, we have $a^{\gamma}e_i\in S$ and thereby $a^{\gamma}\in S_B$. Hence, the desired expression $(*)$ is obtained.}  \par

(ii) By the definition of $S_B$, $S_E$, $L_D$ and the assumption on $S_B$, the property $(**)$ is clear. Thus for a nonzero element $\xi\in N\cap V(S)=N\cap L_D$,  the expression $(*)$ in (i) entails $\xi_1=\xi -\lambda^{-1}\mu a^{\gamma}g\in N\cap L_D$ for some $a^{\gamma}\in S_B\subset D$ and $g\in \G\cap L_D$. If $\xi_1\ne 0$, then repeat this division procedure on $\xi_1$ by $\G\cap L_D$ and so on. As $\prec_e$ is a well-ordering,  after a finite number of repeating the division procedure by $\G\cap M$ we then reach the desired expression for $\xi$. \QED{\parindent=0pt\v5

{\bf Remark} We observe that in the case of Theorem 4.3 (ii), no matter the subalgebra $D$ of $A$ is a solvable polynomial algebra or not (with respect to the restriction of $\prec$ on $S_B$), $\G\cap L_D$ is indeed a \GR basis for the submodule $N\cap L_D$ of $L_D$ (here the restriction of $\prec_e$ on $S$ may be viewed as  a monomial ordering on $S$).}\v5

We now use Theorem 4.3 to derive an analogue of Theorem 2.4 for for submodules of free modules over a solvable polynomial algebra. To see this, Let $A=K[a_1,\ldots ,a_n]$ be a solvable polynomial algebra with admissible system $(\B ,\prec)$, and let $L=\oplus_{i=1}^sAe_i$ be a free $A$-module with admissible system $(\BE, \prec_e)$. For every $1\le d\le s-1$ and the subset $U=\{ e_{i_1},e_{i_2},\ldots ,e_{i_d}\}\subset\{ e_1,\ldots ,e_s\}$ with $i_1<i_2<\cdots <i_d$, let  $L_U=\oplus_{j=1}^{d}Ae_{i_j}$ be the submodule of $L$ generated by $U$. With $S=\{ a^{\alpha}e_{i_j}~|~a^{\alpha}\in\B ,~e_{i_j}\in U\}$ (which is the $K$-basis of $L_U$) and $V(S)=K$-span$S$ (=$L_U$), by Example (1) and Example (2) above we may say that $\prec_e$ is already an elimination ordering on $\BE$ with respect to $V(S)$ (in the sense of Definition 4.1). Moreover, bearing in mind Lemma 4.2, it is straightforward to see that the restriction of $\prec_e$ on $S$ is a monomial ordering on $L_U$.{\parindent=0pt\v5

{\bf 4.4. Theorem}  (Compare with Theorem 2.4.)  With the  preparation made above, if $N$ is a submodule of $L$ and $\G$ is a \GR basis of $N$ with respect to the elimination ordering $\prec_e$, then}\par

(i) $\G\cap L_U$ is a \GR basis for the submodule $N\cap L_U$ of $L_U$ with respect to  the restriction of  $\prec_e$ on $L_U$;\par

(ii) in the case where $N=\sum_{j=1}^mA\xi_j$ is generated by a finite subset  $\{\xi_1,\ldots ,\xi_m\}\subset L$, running the noncommutative Buchberger algorithm for submodules with respect to $\prec_e$ will produce a \GR basis $\G$ for $N$ and consequently, the \GR basis $\G\cap L_U$ for $N\cap L_U$ is obtained.{\parindent=0pt \vskip 6pt

{\bf Proof} (i) With the given subset $S\subset\B (e)$ and the elimination ordering $\prec_e$ on $\BE$, this follows immediately from Theorem 4.3 and the definition of \GR basis for submodules if, in Theorem 4.3, we use $A$ in place of $D$ and use $L_U$ in place of $L_D$.}\par

(ii) This follows from  assertion (i).\QED\v5

By using Theorem 4.3, we also have  a general version of Proposition 2.6 for submodules of free modules (note that an algebra $A$ itself is a free $A$-module). To see this, let $A=K[a_1,\ldots ,a_n]$ be a solvable polynomial algebra with   admissible system $(\B ,\prec )$, and let $L=\oplus_{i=1}^sAe_i$ be a free $A$-module with admissible system $(\BE ,\prec_e)$. Considering the polynomial extension $A[t]$ of $A$ by a commuting variable $t$, let $L_t=\oplus_{i=1}^sA[t]e_i$ be the free $A[t]$-module with the admissible system $(\B (t,e),\prec_{t\textrm{\tiny -}e})$, where
$$\B (t,e)=\{ t^qa^{\alpha}e_i~|~a^{\alpha}\in\B,~q\in\NZ ,~1\le i\le s\} $$ is the $K$-basis of $L_t$ and, with $S=\B (e)$ and thus $V(S)=K$-span$S=L$ (note that $L$ is  an $A$-submodule of $L_t$),  $\prec_{t\textrm{\tiny -}e}$ is the elimination ordering on $\B (t,e)$ with respect to $V(S)$ (as defined in previous Example (4)). Moreover, for any two given submodules $N_1$ and $N_2$ of $L$, let $$N=\sum_{u\in N_1}A[t]tu+\sum_{v\in N_2}A[t](1-t)v$$
be the submodule of $L_t$ generated by $\{ tu, (1-t)v~|~u\in N_1,~v\in N_2\}$. {\parindent=0pt\v5

{\bf 4.5. Proposition}  With the preparation made above, the following statements hold.}\par

(i)  $N\cap L=N_1\cap N_2$.\par

(ii) If $\G$ is a \GR basis of $N$ with respect to the monomial ordering $\prec_{t\textrm{\tiny -}e}$ on $L_t$, then $\G\cap L$ is a \GR basis for $N\cap L=N_1\cap N_2$ with respect to $\prec_e$ on $L$.\par

(iii) In the case where $N_1=\sum^s_{i=1}A\xi_i$ is generated by a finite subset $\{\xi_1,\ldots ,\xi_s\}\subset L$ and $N_2=\sum^m_{j=1}A\eta_j$ is generated by a finite subset $\{\eta_1,\ldots ,\eta_m\}\subset L$, running the noncommutative Buchberger algorithm for submodules with respect to $\prec_{t\textrm{\tiny -}e}$ will produce a \GR basis $\G$ for the submodule  $$N=\sum^s_{i=1}A[t]t\xi_i+\sum_{j=1}^mA[t](1-t)\eta_j$$ of $L_t$ and consequently, $\G\cap L$ gives a \GR basis  for $N\cap L=N_1\cap N_2$ with respect to $\prec_e$ on $L$, which is certainly a generating set of $N_1\cap N_2$.
{\parindent=0pt\vskip 6pt

{\bf Proof} (i) If $\xi\in N_1\cap N_2$,  then $\xi =t\xi+(1-t)\xi\in N$, showing $N_1\cap N_2\subseteq N\cap L$. On the other hand, if $\xi\in N\cap L$, say $\xi=t\theta+(1-t)\eta$ with $\theta\in\sum_{u\in N_1}A[t]u$ and $\eta\in\sum_{v\in N_2}A[t]v$, then, considering the $A$-module homomorphism
$$\begin{array}{cccc} \varphi_1:&L_t&\mapright{}{}&L\\
&\sum_if(t)e_i&\mapsto&\sum_if(0)e_i\end{array}$$
$\xi\in N\cap L$ entails $\xi =\varphi_1(\xi )=\varphi_1(\eta )\in N_2$, and considering the algebra homomorphism
$$\begin{array}{cccc} \varphi_2:&L_t&\mapright{}{}&L\\
&\sum_if(t)e_i&\mapsto&\sum_if(1)e_i\end{array}$$
$\xi\in N\cap L$ entails $\xi =\varphi_2(\xi )=\varphi_2(\theta )\in N_1$. Hence $\xi\in N_1\cap N_2$ and consequently $N\cap L\subseteq N_1\cap N_2$. Combining both inclusions we conclude  $N_1\cap N_2=N\cap L$.}\par

(ii) With the given  subset $S=\B (e)\subset\B (t,e)$ and the elimination ordering  $\prec_{t\textrm{\tiny -}e}$ on $\B (e,t)$, this follows immediately from Theorem 4.3 and the definition of \GR basis for submodules if, in Theorem 4.3, we use $A$ in place of $D$ and use $L$ in place of $L_D$.\par

(iii) This follows from  assertion (ii).\v5

\section*{5. Applications to Module Homomorphisms}
Let $A=K[a_1,\ldots ,a_n]$ be a solvable polynomial algebra with admissible system $(\B ,\prec )$. In this section we apply the results of previous sections to $A$-module homomorphisms in order to{\parindent=1.3truecm\par

\item{(1)} get a generating set for the kernel of $A$-module homomorphisms between free $A$-modules, respectively for the kernel of $A$-module homomorphisms between quotient modules of free modules (thus the injectivity of the homomorphisms concerned may be determined);\par

\item{(2)} solve the membership problem for the image of  $A$-module homomorphisms  between free $A$-modules, respectively for the image of $A$-module homomorphisms between quotient modules of free modules (thus the surjectivity of the homomorphisms concerned may be determined).\par} \v5

Let $L_1\oplus_{i=1}^sAe_i$ and $L_2=\oplus_{j=1}^mA\VE_j$ be free $A$-modules, and let $\varphi$: $L_1\r L_2$ be an $A$-module homomorphism such that $\varphi (e_i)=\eta_i\in L_2$, $1\le i\le s$. Our first goal is to get a generating set for Ker$\varphi$ (the kernel of $\varphi$) by computing a \GR basis for it, and to solve the membership problem for the image Im$\varphi$ of $\varphi$.
{\parindent=0pt\v5

{\bf 5.1. Lemma}   With the $A$-module homomorphism $\varphi$ defined above, consider the free $A$-module $L=L_1\oplus L_2$ and the $A$-module homomorphism $\Phi$: $L\r L_2$ such that $\Phi (\VE_j)=\VE_j$, $1\le j\le m$, and $\Phi (e_i)=\eta_i$, $1\le i\le s$. Then }

(i) the diagram of $A$-module homomorphisms
$$\begin{diagram} L_1&\rTo^{\scriptstyle\varphi}&L_2\\
\scriptstyle{\iota~}\dTo&\NE_{\scriptstyle\Phi}\\
L&&\end{diagram}$$
is commutative, i.e., $\Phi\circ\iota =\varphi$, where $\iota$ is the inclusion map;\par

(ii) $\hbox{Ker}\Phi =\sum^s_{i=1}A(e_i-\eta _i).${\parindent=0pt\vskip 6pt

{\bf Proof}  (i) By the definition of each homomorphism in the diagram, it is clear that the diagram is commutative.}\par

(ii) First note that the definition of $\Phi$ entails that $\Phi (\xi )=\varphi (\xi )$ for all $\xi \in L_1$, $\Phi (\eta )=\eta$ for all $\eta\in L_2$. It follows that $\sum^s_{i=1}A(e_i-\eta _i)\subseteq$ Ker$\Phi$. On the otherf hand, if $x=\xi+\eta\in$ Ker$\Phi$, where $\xi=\sum_if_ie_i\in L_1$ and $\eta\in L_2$, then
$$\renewcommand\arraystretch{1.5}\begin{array}{rcl} 0=\Phi (x)&=&\Phi \left (\displaystyle\sum_if_ie_i+\eta \right )\\
&=&\displaystyle\sum_if_i\Phi (e_i)+\Phi (\eta )\\
 &=&\displaystyle\sum_if_i\varphi (e_i)+\eta \\
&=&\displaystyle\sum_if_i\eta_i+\eta.\end{array}$$
It turns out that
$$\renewcommand\arraystretch{1.5}\begin{array}{rcl} x&=&\xi +\eta\\
&=&\displaystyle\sum_if_ie_i-\sum_if_i\eta_i\\
&=&\displaystyle\sum_if_i(e_i-\eta_i)\in \displaystyle\sum^s_{i=1}A(e_i-\eta _i),\end{array}$$
and thus Ker$\Phi\subseteq \sum^s_{i=1}A(e_i-\eta _i)$. Combining the inclusions of both directions, we conclude that Ker$\Phi = \sum^s_{i=1}A(e_i-\eta _i)$.\QED {\parindent=0pt\v5

{\bf 5.2. Theorem}  Let the data $(L_1, L_2, \varphi ; L,\Phi , \hbox{Ker}\Phi )$ be as in Lemma 5.1. Then Ker$\varphi =$ Ker$\Phi\cap L_1= (\sum^s_{i=1}A(e_i-\eta _i))\cap L_1$. \vskip 6pt

{\bf Proof} By the definitions of $\varphi$ and $\Phi$ we have $\Phi (\xi )=\varphi (\xi )$ for all $\xi\in L_1$, thereby Ker$\varphi\subseteq$ Ker$\Phi\cap L_1$. By Lemma 5.1(i) we know that $\Phi\circ\iota =\varphi$. So, for $x\in$ Ker$\Phi\cap L_1$ we have $0=\Phi (x)=(\Phi\circ\iota )(x)=\varphi (x)$, thus Ker$\Phi\cap L_1\subseteq$ Ker$\varphi$. Therefore, Ker$\varphi =$ Ker$\Phi\cap L_1= (\sum^s_{i=1}A(e_i-\eta _i))\cap L_1$.\QED}\v5

For the convenience of our next usage, let us rewrite $L_2$ as $L_2=\oplus_{j=s+1}^{m+s}A\VE_j$ so that  $L=L_1\oplus L_2=(\oplus_{i=1}^sAe_i)\oplus (\oplus_{j=s+1}^{m+s}A\VE_j)$. {\parindent=0pt\v5

{\bf 5.3. Corollary}   With the data $(L_1, L_2, \varphi ,\hbox{Ker}\VF ; L,\Phi , \hbox{Ker}\Phi )$  as in Lemma 5.1, Theorem 5.2 and that fixed above,  let $\prec_{pot}$ be  the POT  elimination monomial ordering on the $K$-basis $\B (e,\VE )=\B (e)\cup\B (\VE )$ of $L$ with respect to $V(S)$ where $S=\B (e)$ is the $K$-basis of $L_1$ and $V(S)=K$-span$S=L_1$ (see Example (1) of Section 4), and let $\G$ be a \GR basis of  Ker$\Phi =\sum^s_{i=1}A(e_i-\eta _i)$ in $L$ computed by running  the noncommutative Buchberger algorithm for submodules with respect to $\prec_{pot}$, then $\G\cap L_1$ is a \GR basis of Ker$\varphi =$ Ker$\Phi\cap L_1$ in $L_1$ with respect to the restriction of $\prec_{pot}$ on $L_1$. \vskip 6pt

{\bf Proof} This follows from Theorem 4.4.\QED}\v5

We next turn to solve the membership problem for the image Im$\varphi$ of the $A$-module homomorphisms $L_1~\mapright{\varphi}{}~L_2$, where $\varphi (e_i)=\eta_i$ for $1\le i\le s$. To this end, we keep using the notation fixed above, in particular, $L=L_1\oplus L_2=(\oplus_{i=1}^sAe_i)\oplus (\oplus_{j=s+1}^{m+s}A\VE_j)$ with the $K$-basis $\B (e,\VE )=\B (e)\cup\B (\VE )$, $L~\mapright{\Phi}{}~L_2$ with $\Phi (e_i)=\eta_i=\varphi (e_i)$, $1\le i\le s$, $\Phi (\VE_j)=\VE_j$ with $s+1\le j\le m+s$, and Ker$\Phi = \sum^s_{i=1}A(e_i-\eta _i)$. {\parindent=0pt\v5

{\bf 5.4. Theorem}  With the notation above, let $\prec_{pot}$ be the POT elimination ordering on $\B (e,\VE )$ with respect to $V(S)$, where $S=\B (e)$ and $V(S)=K$-span$S=L_1$ (see Example (1) of Section 4),  and let $\G$ be a \GR basis of Ker$\Phi$ in $L$ with respect to  $\prec_{pot}$.  For an element $\eta\in L_2$, that $\eta\in $ Im$\varphi$ if and only if there is a $\xi =\sum_if_ie_i\in L_1$ such that $\OV{\eta}^{\G}=\xi$. If this is the case, then $\eta =\varphi (\xi )=\sum_if_i\eta_i$.\vskip 6pt

{\bf Proof} If $\eta\in$ Im$\VF$, then there is a $\xi^*\in L_1$ such that $\VF (\xi^*)=\eta$. Thus, by Lemma 5.1 we have $\Phi (\xi^*)=\VF (\xi^*)=\eta =\Phi (\eta )$ and thereby $\eta -\xi^*\in$ Ker$\Phi$. Since $\G$ is a \GR basis of Ker$\Phi$ with respect to $\prec_{pot}$, it follows that $0=\OV{\eta -\xi^*}^{\G}=\OV{\eta}^{\G}-\OV{\xi^*}^{\G}$. If $\OV{\xi^*}^{\G}=0$, then $\OV{\eta}^{\G}=0$. If $\OV{\xi^*}^{\G}\ne 0$, then noticing that $\xi^*\in L_1$ and $\prec_{pot}$ is an elimination ordering with respect to $V(S)$ where $S=\BE$ and $V(S)=K$-span$S=L_1$, the property $\LM (\OV{\xi^*}^{\G})\prec_{pot}\LM (\xi^*)$ of a remainder entails $\OV{\xi^*}^{\G}\in L_1$. Putting $\xi =\OV{\xi^*}^{\G}$, we then have $\OV{\eta}^{\G}=\xi\in L_1$.}\par

Conversely, if there is a $\xi=\sum_if_ie_i\in L_1$ such that $\OV{\eta}^{\G}=\xi$, then since a division procedure on $\eta$ by $\G$ yields $\eta =\sum_qh_qg_q+\OV{\eta}^{\G}$, where $h_q\in A$, $g_q\in\G\subset$ Ker$\Phi$, it follows from Lemma 5.1 that
$$\begin{array}{rcl} \eta =\Phi (\eta )&=&\sum_qh_q\Phi (g_q)+\Phi (\OV{\eta}^{\G})\\
&=&\Phi (\xi )=\VF (\xi )=\sum_if_i\eta_i\in~\hbox{Im}\VF ,\end{array}$$
finishing the proof.\QED{\parindent=0pt\v5

{\bf 5.5. Corollary}  With the notation as used in Theorem 5.4, for an element $\eta\in L_2$, that $\eta\in$ Im$\VF$ if and only if $\OV{\eta}^{\G}\in L_1$.\par\QED\v5

{\bf 5.6. Theorem}  With the notation as used in Theorem 5.4, let $\G$ be a reduced \GR basis of Ker$\Phi$ in $L$ with respect to the elimination ordering $\prec_{pot}$ on $L$. Then $\VF$ is surjective, that is, Im$\VF =\sum^s_{i=1}A\eta_i=\sum^{m+s}_{j=s+1}A\VE_j=L_2$, if and only if for each $j=s+1,\ldots , m+s$, there is a $g_j\in\G$ such that $g_j=\VE_j-\xi_j$ for some $\xi_j=\sum^s_{i=1}f_{ji}e_i\in L_1$. If this is the case, then $\VE_j=\sum^s_{i=1}f_{ji}\eta_i$ for $s+1\le j\le m+s$.\vskip 6pt

{\bf Proof} If $\VF$ is surjective, then every $\VE_j\in$ Im$\VF$, $s+1\le j\le m+s$. It follows from Theorem 5.4 that for each $j=s+1,\ldots , m+s$, there is a $\xi_j'\in L_1$ such that $\OV{\VE_j}^{\G}=\xi_j'$. In other words, after implementing the division procedure on $\VE_j$ by $\G$ we have $\VE_j=\sum_th_tg_i+\OV{\VE_j}^{\G}=\sum_th_tg_t+\xi_j'$, where $h_t\in A$ and $g_t\in\G$. Since every $\xi_j'\in L_1$ while $\prec_{pot}$ is the elimination ordering with respect to $V(S)=K$-span$\B (e)=L_1$, it follows that $\LM (\xi_j')\prec_{pot}\VE_j$ and thus there is some $t$ such that $\VE_j=\LM (\VE_j)=\LM (h_tg_t)$. But this implies $\LM (g_t)|_{_{\rm L}}\VE_j$, thereby $\LM (g_t)=\VE_j$ (see Lemma 3.2). Without loss of generality we may say that $t=j$. Hence  $\LM (g_j)=\VE_j$, $s+1\le j\le m+s$. Furthermore, noticing that $\G$ is a reduced \GR basis for Ker$\Phi$, there must exist $\xi_j\in L_1$ such that $g_j=\VE_j-\xi_j$ for every $s+1\le j\le m+s$.}\par

Conversely, if for each $j=s+1,\ldots , m+s$, there is a $g_j\in\G\subset$ Ker$\Phi$ such that $g_j=\VE_j-\xi_j$ for some $\xi_j=\sum^s_{i=1}f_{ji}e_i\in L_1$, then
$$\begin{array}{rcl} 0=\Phi (g_j)&=&\Phi (\VE_j)-\Phi (\xi_j)\\
&=&\VE_j-\sum_if_{ji}\Phi (e_i)\\
&=&\VE_j-\sum^s_{i=1}f_{ji}\VF (e_i)\\
&=&\VE_j-\sum^s_{i=1}f_{ji}\eta_i\end{array}$$
and consequently $\VE_j=\sum^s_{i=1}f_{ji}\eta_i\in$ Im$\VF$. This shows that $\VF$ is surjective.\par\QED\v5

Let $A=K[a_1,\ldots ,a_n]$ be a solvable polynomial algebra with admissible system $(\B ,\prec )$. Now, as dealing with $A$-module homomorphisms between free modules, we proceed to derive similar results for $A$-module homomorphisms between quotient modules of free modules.\par

Rcall that if  $M_1=\sum^s_{i=1}A\theta_i$ is an $A$-module generated by a finite subset $\{\theta_1,\ldots ,\theta_s\}\subset M$ and $L_1=\oplus_{i=1}^s$ is the free $A$-module with  $A$-basis $\{ e_1,\ldots ,e_s\}$, then there is an $A$-module epimorphism $\psi$: $L_1\r M$ such that $\psi (e_i)=\theta_i$, $1\le i\le s$, and thus $M_1\cong L_1/N_1$ with $N_1=$ Ker$\psi$. This shows that every finitely generated $A$-module can be presented as a quotient module of a free module. Let $M_2=\sum_{j=s+1}^{m+s}A\nu_j$ be another  $A$-module generated by a finite subset $\{\nu_{s+1},\ldots ,\nu_{m+s}\}\subset M_2$, then similarly $M_2\cong L_1/N_2$, where $L_2=\oplus_{j=s+1}^{m+s}A\VE_j$ is the free $A$-module with $A$-basis $\{\VE_{s+1},\ldots ,\VE_{m+s}\}$ and $N_2$ is a submodule of $L_2$. We first have an criterion for the existence of an $A$-module homomorphism from $M_1$ to $M_2$.{\parindent=0pt\v5

{\bf 5.7. Proposition}  With the notation fixed above, if $N_1=\sum^t_{q=1}A\xi_q$ with each $\xi_q=\sum^s_{i=1}f_{qi}e_i\in N_1$, and $N_2=\sum_{\ell =1}^dA\omega_{\ell}$ with every $\omega_{\ell}\in N_2$, then there exists an $A$-module homomorphism $\VF$: $M_1\r M_2$ if and only if there are $\eta_1,\ldots ,\eta_s\in L_2$ such that $\sum^s_{i=1}f_{qi}\eta_i\in N_2$ holds for $1\le q\le t$ (note that this membership problem can be solved by using a \GR basis of $N_2$).\vskip 6pt

{\bf Proof} Note that for any given  $\eta_1,\ldots ,\eta_s\in L_2$, there exists an $A$-module homomorphism $\phi$: $L_1\r L_2$ such that $\phi (e_i)=\eta_i$ for $1\le i\le s$. So, considering  the canonical homomorphism $\pi_1$: $L_1\r M_1$, the canonical homomorphism $\pi_2$: $L_2\r M_2$, and  the following diagram of $A$-module homomorphisms:
$$\begin{diagram} L_1&\rTo^{\scriptstyle\phi}{}&L_2\\
\scriptstyle{\pi_1~}\dTo&&\dTo{}{\scriptstyle \pi_2}\\
M_1=L_1/N_1&\rTo_{\scriptstyle ?}&M_2=L_2/N_2\end{diagram}$$
one may check that  the given condition is a sufficient and necessary condition for having an $A$-module homomorphism $\OV{\phi}$: $M_1\r M_2$ such that the diagram is commutative.\QED}\v5

Before continuing, let us make a convention for convenience, namely if $L=\oplus_{i=1}^sAe_i$ is a free $A$-module and $N$ is a submodule of $L$, then for $\xi\in L$ we write $\OV{\xi}$ for the coset in $L/N$ represented by $\xi$.{\parindent=0pt\v5

{\bf 5.8. Lemma}  (Compare with Lemma 5.1.) Let $N_1$ be a submodule of the free $A$-module  $L_1=\oplus_{i=1}^sAe_i$, $N_2$ a submodule of the free $A$-module $L_2=\oplus_{j=s+1}^{m+s}A\VE_j$, and let $L=L_1\oplus L_2$. Considering the $A$-modules $M_1=L_1/N_1$, $M_2=L_2/N_2$ and an $A$-module homomorphism $\varphi$: $M_1\r M_2$ such that $\VF (\OV{e_i})=\OV{\eta_i}$ with each  $\eta_i\in L_2$, let us take the canonical homomorphism $\pi_e$: $L_1\r M_1=L_1/N_1$ and the $A$-module homomorphism $\Phi$: $L\r M_2$ such that $\Phi (e_i)=\OV{\eta_i}=\varphi (\OV{e_i})$, $1\le i\le s$, and $\Phi (\VE_j)=\OV{\VE_j}$, $s+1\le j\le m+s$. Then}\par

(i) the diagram of $A$-homomorphisms
$$\begin{diagram} &L&\lTo^{\scriptstyle\iota}&L_1\\
&{\scriptstyle\Phi~}\dTo&&~\dTo_{\scriptstyle \pi_e}\\
L_2/N_2=&M_2&\lTo{\scriptstyle\VF}&M_1&=L_1/N_1\end{diagram}$$
is commutative, i.e., $\VF\circ\pi_e=\Phi\circ\iota$, where $\iota$ is the inclusion map.\par

(ii) Ker$\Phi =N_2+\sum_{i=1}^sA(e_i-\eta_i)$.{\parindent=0pt\vskip 6pt

{\bf Proof} A direct verification is left as an exercise.\QED\v5

{\bf 5.9. Theorem}  (Compare with Theorem 5.2.) With the data $(M_1,M_2,\VF ;L,\Phi, \hbox{Ker}\Phi)$ as in Lemma 5.8,
$$\hbox{Ker}\VF =(\hbox{Ker}\Phi\cap L_1+N_1)/N_1.$$\v5

{\bf Proof} A direct verification is left as an exercise.\QED\v5

{\bf 5.10. Corollary}   With the data $(M_1,M_2,\VF , \hbox{Ker}\VF ;L,\Phi, \hbox{Ker}\Phi)$ as in Lemma 5.8 and Theorem 5.9, let $\prec_{pot}$ be  the POT  elimination monomial ordering on the $K$-basis $\B (e,\VE )=\B (e)\cup\B (\VE )$ of $L$ with respect to $V(S)$ where $S=\B (e)$ is the $K$-basis of $L_1$ and $V(S)=K$-span$S=L_1$ (see Example (1) of Section 4). If   $N_2=\sum_{\ell =1}^dA\omega_{\ell}$ with every $\omega_{\ell}\in N_2$, and $\G$ is a \GR basis of  Ker$\Phi =N_2+\sum^s_{i=1}A(e_i-\eta _i)$ in $L$ computed by running  the noncommutative Buchberger's algorithm with respect to $\prec_{pot}$, then $\G\cap L_1$ gives rise to a generating set of Ker$\varphi =($ Ker$\Phi\cap L_1+N_1)/N_1$ in $M_1$. \vskip 6pt

{\bf Proof} This follows from Theorem 4.4.\QED}\v5

We next turn to solve the membership problem for the image Im$\varphi$ of the $A$-module homomorphisms $M_1=L_1/N_1~\mapright{\varphi}{}~M_2=L_2/N_2$, where $\varphi (\OV{e_i})=\OV{\eta_i}$ for $1\le i\le s$. To this end, we keep using the notation fixed above, in particular, $L=L_1\oplus L_2=(\oplus_{i=1}^sAe_i)\oplus (\oplus_{j=s+1}^{m+s}A\VE_j)$ with the $K$-basis $\B (e,\VE )=\B (e)\cup\B (\VE )$, $L~\mapright{\Phi}{}~L_2/N_2$ with $\Phi (e_i)=\OV{\eta_i}=\varphi (\OV{e_i})$, $1\le i\le s$, $\Phi (\VE_j)=\OV{\VE_j}$ with $s+1\le j\le m+s$, and Ker$\Phi = N_2+\sum^s_{i=1}A(e_i-\eta _i)$. {\parindent=0pt\v5

{\bf 5.11. Theorem}  (Compare with Theorem 5.4.) With the notation above, let $\prec_{pot}$ be the POT elimination ordering on $\B (e,\VE )$ with respect to $V(S)$, where $S=\B (e)$ and $V(S)=K$-span$S=L_1$ (see Example (1) of Section 4),  and let $\G$ be a \GR basis of Ker$\Phi$ in $L$ with respect to  $\prec_{pot}$.  For an element $\OV{\eta}\in M_2$, that $\OV{\eta}\in $ Im$\varphi$ if and only if there is a $\xi =\sum_if_ie_i\in L_1$ such that $\OV{\eta}^{\G}=\xi$. If this is the case, then $\OV{\eta} =\varphi (\OV{\xi} )=\sum_if_i\OV{\eta_i}$.\vskip 6pt

{\bf Proof} If $\OV{\eta}\in$ Im$\VF$, then there is a $\OV{\xi^*}\in M_1=L_1/N_1$ such that $\VF (\OV{\xi^*})=\OV{\eta}$. Thus, by Lemma 5.8 we have $\Phi (\xi^*)=\VF (\OV{\xi^*})=\OV{\eta} =\Phi (\eta )$ and thereby $\eta -\xi^*\in$ Ker$\Phi\subset L$. Since $\G$ is a \GR basis of Ker$\Phi$ with respect to $\prec_{pot}$, it follows that $0=\OV{\eta -\xi^*}^{\G}=\OV{\eta}^{\G}-\OV{\xi^*}^{\G}$. If $\OV{\xi^*}^{\G}=0$, then $\OV{\eta}^{\G}=0$. If $\OV{\xi^*}^{\G}\ne 0$, then noticing that $\xi^*\in L_1$ and $\prec_{pot}$ is an elimination ordering with respect to $V(S)$ where $S=\BE$ and $V(S)=K$-span$S=L_1$, the property $\LM (\OV{\xi^*}^{\G})\prec_{pot}\LM (\xi^*)$ of a remainder entails $\OV{\xi^*}^{\G}\in L_1$. Putting $\xi =\OV{\xi^*}^{\G}$, we then have $\OV{\eta}^{\G}=\xi\in L_1$.}\par

Conversely, if there is a $\xi=\sum_if_ie_i\in L_1$ such that $\OV{\eta}^{\G}=\xi$, then since a division procedure on $\eta$ by $\G$ yields $\eta =\sum_qh_qg_q+\OV{\eta}^{\G}$, where $h_q\in A$, $g_q\in\G\subset$ Ker$\Phi$, it follows from Lemma 5.8 that
$$\begin{array}{rcl} \OV{\eta} =\Phi (\eta )&=&\sum_qh_q\Phi (g_q)+\Phi (\OV{\eta}^{\G})\\
&=&\Phi (\xi )=\VF (\OV{\xi} )=\sum_if_i\OV{\eta_i}\in~\hbox{Im}\VF ,\end{array}$$
finishing the proof.\QED{\parindent=0pt\v5

{\bf 5.12. Corollary}  With the notation as used in Theorem 5.11, for an element $\OV{\eta}\in M_2$, that $\OV{\eta}\in$ Im$\VF$ if and only if $\OV{\eta}^{\G}\in L_1$.\par\QED\v5

{\bf 5.13. Theorem}  (Compare with Theorem 5.6.) With the notation as used in Theorem 5.11, let $\G$ be a reduced \GR basis of Ker$\Phi$ in $L$ with respect to the elimination ordering $\prec_{pot}$ on $L$. Then $\VF$ is surjective, that is, Im$\VF =\sum^s_{i=1}A\OV{\eta_i}=\sum^{m+s}_{j=s+1}A\OV{\VE_j}=M_2$, if and only if for each $j=s+1,\ldots , m+s$, there is a $g_j\in\G$ such that $g_j=\VE_j-\xi_j$ for some $\xi_j=\sum^s_{i=1}f_{ji}e_i\in L_1$. If this is the case, then $\OV{\VE_j}=\sum^s_{i=1}f_{ji}\OV{\eta_i}$ for $s+1\le j\le m+s$.\vskip 6pt

{\bf Proof} If $\VF$ is surjective, then every $\OV{\VE_j}\in$ Im$\VF$, $s+1\le j\le m+s$. It follows from Theorem 5.11 that for each $j=s+1,\ldots , m+s$, there is a $\xi_j'\in L_1$ such that $\OV{\VE_j}^{\G}=\xi_j'$. In other words, after implementing the division procedure on $\VE_j$ by $\G$ we have $\VE_j=\sum_th_tg_i+\OV{\VE_j}^{\G}=\sum_th_tg_t+\xi_j'$, where $h_t\in A$ and $g_t\in\G$. Since every $\xi_j'\in L_1$ while $\prec_{pot}$ is the elimination ordering with respect to $V(S)=K$-span$\B(e)=L_1$, it follows that $\LM (\xi_j')\prec_{pot}\VE_j$ and thus there is some $t$ such that $\VE_j=\LM (\VE_j)=\LM (h_tg_t)$. But this implies $\LM (g_t)|_{_{\rm L}}\VE_j$, thereby $\LM (g_t)=\VE_j$ (see Lemma 3.2). Without loss of generality we may say that $t=j$ and hence $\LM (g_j)=\VE_j$, $s+1\le j\le m+s$. Furthermore, noticing that $\G$ is a reduced \GR basis for Ker$\Phi$, there must exist $\xi_j\in L_1$ such that $g_j=\VE_j-\xi_j$ for every $s+1\le j\le m+s$.}\par

Conversely, if for each $j=s+1,\ldots , m+s$, there is a $g_j\in\G\subset$ Ker$\Phi$ such that $g_j=\VE_j-\xi_j$ for some $\xi_j=\sum^s_{i=1}f_{ji}e_i\in L_1$, then
$$\begin{array}{rcl} 0=\Phi (g_j)&=&\Phi (\VE_j)-\Phi (\xi_j)\\
&=&\OV{\VE_j}-\sum_if_{ji}\Phi (e_i)\\
&=&\OV{\VE_j}-\sum^s_{i=1}f_{ji}\VF (\OV{e_i})\\
&=&\OV{\VE_j}-\sum^s_{i=1}f_{ji}\OV{\eta_i}\end{array}$$
and consequently $\OV{\VE_j}=\sum^s_{i=1}f_{ji}\OV{\eta_i}\in$ Im$\VF$. This shows that $\VF$ is surjective.\par\QED\v5

\centerline{\bf References} {\parindent=.8truecm\vskip 8pt

\item{[AL]} W. W. Adams and P. Loustaunau, {\it An Introduction to Gr\"obner Bases}.
Graduate Studies in Mathematics, Vol. 3. American Mathematical
Society, 1994.

\item{[Ch]} F.~Chyzak, Holonomic systems and automatic proving of
identities, {\it Research Report} 2371,  Institute National de
Recherche en Informatique et en Automatique, 1994.

\item{[CS]} F. Chyzak and B. Salvy, Noncommutative elimination in Ore
algebras proves multivariate identities, {\it J. Symbolic Comput}.,
26(1998), 187--227.

\item{[Eis]} D. Eisenbud, {\it Commutative Algebra with a View
Toward to Algebraic Geometry}, GTM 150. Springer, New York, 1995.

\item{[JBS]} A. Jannussis, et al,  Remarks on
the $q$-quantization. {\it Lett. Nuovo Cimento}, 30(1981), 123--127.

\item{[KR]} M. Kreuzer, L. Robbiano, {\it Computational Commutative Algebra 1}. Springer, 2000.

\item{[K-RW]} A.~Kandri-Rody and V.~Weispfenning, Non-commutative
Gr\"obner bases in algebras of solvable type. {\it J. Symbolic
Comput.}, 9(1990), 1--26.

\item{[KW]} H. Kredel and V. Weispfenning, Computing dimension and independent sets for polynomial Ideals. {\it J. Symb. Comput}. 6(2/3)(1988), 231--247.

\item{[Kur]} M.V. Kuryshkin, Op\'erateurs quantiques
g\'en\'eralis\'es de cr\'eation et d'annihilation. {\it Ann. Fond.
L. de Broglie}, 5(1980), 111--125.

\item{[Lev]} V. Levandovskyy, {\it Non-commutative Computer Algebra for
Polynomial Algebra}: {\it Gr\"obner Bases, Applications and
Implementation}. Ph.D. Thesis, TU Kaiserslautern, 2005.

\item{[Li1]} H. Li, {\it Noncommutative Gr\"obner Bases and
Filtered-Graded Transfer}. Lecture Note in Mathematics, Vol. 1795,
Springer, 2002.

\item{[Li2]} H. Li, {\it Gr\"obner Bases in Ring Theory}. World Scientific
Publishing Co., 2011.

\item{[Li3]} H. Li, A note on solvable polynomial algebras.
{\it Computer Science Journal of Moldova}, 1(64)(2014), 99--109.

\item{[Li4]} H. Li, An elimination lemma for algebras with PBW bases. {\it Communications in Algebra}, 8(46)(2018), 3520¨C3532.

\item{[Li5]} H. Li, Notes on \GR bases and free resolutions of modules over solvable polynomial algebras.
arXiv:1510.04381v1[math.RA], 131 pages.

\item{[LW]} H. Li and  Y. Wu, ~Filtered-graded transfer of
Gr\"obner basis computation in solvable polynomial algebras. {\it
Communications in Algebra}, 1(28)(2000), 15--32.

\item{[PWZ]} M. Petkov¡¦sek, H. Wilf and D. Zeilberger, {\it $A = B  $}. A.K.
Peters, Ltd. 1996.

\end{document}